\documentclass[10pt]{article}
\usepackage{amsmath,amsfonts,amssymb}
\usepackage{graphics} 
\usepackage{epsfig} 
\usepackage{psfrag}

\newtheorem{lemma}{Lemma}
\newtheorem{theorem}{Theorem}

\newtheorem{corollary}{Corollary}
\newtheorem{defn}{Definition}

\newtheorem{remark}{Remark}

\newtheorem{conjecture}{Conjecture}
\newtheorem{algorithm}{Algorithm}

\author{Mark Verwoerd and Oliver Mason\\\\Hamilton Institute, National University of Ireland\\\\{\tt \{mark.verwoerd,oliver.mason\}@nuim.ie}}
\title{Global phase-locking in finite populations of
phase-coupled oscillators}
\begin{document}
\maketitle

\def\QEDclosed{\mbox{\rule[0pt]{1.3ex}{1.3ex}}} 
\def\QEDopen{{\setlength{\fboxsep}{0pt}\setlength{\fboxrule}{0.2pt}\fbox{\rule[0pt]{0pt}{1.3ex}\rule[0pt]{1.3ex}{0pt}}}}
\def\QED{\QEDclosed} 
\def\proof{\noindent\hspace{2em}{\itshape Proof: }}
\def\endproof{\hspace*{\fill}~\QED\par\endtrivlist\unskip}

\begin{abstract}

We present new necessary and sufficient conditions for the
existence of fixed points in a finite system of coupled
phase oscillators on a complete graph. We use these
conditions to derive bounds on the critical coupling.
\end{abstract}

\section{Introduction}
The phenomenon of synchronization arises in a wide variety
of application areas across neuroscience, biology,
engineering and physics \cite{GolSte06, Pik03, Beu03,
Mur02, Gla01}.  As such, the identification and study of
structures and mechanisms that support the onset of
synchronized behaviour is a key issue in the theory of
interconnected dynamical systems.  In particular, there has
been a great deal of interest across the mathematics,
physics and engineering communities in the development and
analysis of simple mathematical models of synchronization
\cite{Jad04,Rog04a,Rog04b,Cho05,Sep06,SPL06}.

To date, one of the most widely-studied frameworks for the
analysis of synchronization is the so-called {\it Kuramoto
model} of phase coupled oscillators \cite{Kur75, Stro00}.
In fact, this model has been used in numerous applications
in the chemical and biological sciences, and its basic
properties have been analysed using a combination of
numerical and analytical techniques
\cite{Kur84,Stro00,Stro01,Ace05}.  The basic Kuramoto model
is comprised of a system of coupled oscillators, which may
have different natural frequencies, where the coupling
between two oscillators is given by a weighted sinusoidal
function of the difference of their phases.  The weights
used in the model are typically taken to be the same for
all pairs of oscillators and are given by the ratio of a
fixed parameter, the coupling strength, to the network
size.

The aspect of the Kuramoto model that has attracted most
attention to date is the manner in which the onset of
synchronization depends on the strength of coupling between
the oscillators. For instance, at very low values of the
coupling strength, little or no synchronization is
observed. As the coupling strength is increased, some
partial synchronization appears in the network up to a
threshold value of the coupling strength, referred to here
as the critical coupling, at which fully synchronized
behaviour emerges \cite{Jad04,Cho05}. The mechanism of
(de)synchronization in finite populations of oscillators
has been described in considerable detail
\cite{MaiPop04,MaiPop05b}. In particular, when the coupling
strength drops below its critical value and as it continues
to decrease, the system undergoes a series of so called
\emph{frequency-splitting bifurcations}. At each such
bifurcation, the ensemble of oscillators subdivides into
smaller and smaller groups of oscillators with identical
average frequency, until eventually all oscillators
oscillate at their own intrinsic frequency. A detailed
analysis of this behaviour for a system with three
oscillators was given in \cite{MaiPop04}. While the
aforementioned contributions focus on the behaviour of the
system in the \emph{sub}critical coupling regime, the
present paper studies globally phase-locked solutions,
which by definition only exist in the \emph{super}critical
coupling regime.

Another aspect of the Kuramoto model to have attracted
attention recently is the emergence of phase
chaos~\cite{MaiPop05,PopMai05} in systems of dimension four
and higher. A generic feature of coupled oscillator
systems, phase chaos in the Kuramoto model is most
prominent in systems with relatively low dimension
(comprising between ten and fifteen
oscillators)~\cite{MaiPop05}. Again this phenomenon can
only exist in the subcritical coupling regime, and we shall
not further consider it here.

 In the original
Kuramoto model, it is assumed that all pairs of oscillators
in the network are connected with the same coupling
strength \cite{Kur75}. This type of coupling is referred to
as `all-to-all' coupling and corresponds to a network in
which the underlying graph is complete \cite{Die00}.
Extensions of the Kuramoto model to lattices \cite{Hong05}
and rings \cite{Rog04b} have also been considered, and more
recently the dynamics of coupled oscillators on networks
with small-world \cite{Watt98,Hong02} and scale-free
\cite{Mor04} topologies have started to attract a lot of
interest.  More generally, there are many fundamental
questions relating to the interplay between a network's
topology and dynamical processes taking place on it which
are still unanswered. The work described in
\cite{GolSte06}, which proposes an extension of group-based
symmetry, using the so-called groupoid formalism, as a
means of classifying possible behaviours for networked
dynamical systems, is particularly noteworthy in this
context.

Many of the recent results concerned with the dynamics and
synchronization of coupled oscillators have either been
based on numerical simulations or else have been derived
for the limiting case of networks of infinite size.  In
contrast, relatively few rigorous results are available for
finite-size networks \cite{Stro00,Jad04}.  In this paper,
we shall be concerned with synchronization in finite
systems of coupled oscillators. Specifically:~we shall
establish (new) necessary and sufficient conditions for the
existence of fixed points in a finite system of coupled
oscillators (see also \cite{Ver06,Ver07}); compute bounds
on the critical coupling strength for such systems; and
provide insights into the number of fixed points possible
under strong coupling. Our analysis is in the spirit of the
work presented in \cite{Jad04,Cho05}, and places particular
emphasis on the \emph{existence} of fixed points.  Of
course, the stability of such fixed points is also a topic
of great interest, and has been considered in
\cite{Jad04,Rog04a,Rog04b,Cho05}.  However, we shall not
explicitly address the question of stability in the current
paper.

The outline of the paper is as follows. In Section~2, we
introduce the Kuramoto model, and review some of its basic
properties. Here, we also give a formal definition of the
critical coupling, which is essentially the lowest value of
the coupling strength for which fixed points exist. In
Section~3, we show that fixed points will always exist for
sufficiently strong coupling (essentially proving that the
critical coupling is a finite number), and then, in
Section~4 provide lower bounds on the critical coupling.
Section~5 contains necessary and sufficient conditions for
the existence of fixed points, which are then used in
Section~6 to describe an algorithm for computing the
critical coupling. Section~7 contains a numerical example
to illustrate the results of the paper and finally, in
Section~8 we present our concluding remarks.

\section{Mathematical Preliminaries and The Kuramoto model}
\subsection{Basic Notation}
Throughout the paper, $\mathbb{R}$ ($\mathbb{C}$) denotes
the field of real (complex) numbers, $\mathbb{R}^N$
($\mathbb{C}^N$) denotes the vector space of all $N$-tuples
of real (complex) numbers, and $\mathbb{R}^{N \times N}$
($\mathbb{C}^{N \times N}$) denotes the space of $N \times
N$ matrices with entries in $\mathbb{R}$ ($\mathbb{C}$).
$i$ is used to denote the complex number satisfying $i^2 =
-1$.  For a vector $x \in \mathbb{R}^N$, $x_i$ denotes the
$i^\text{th}$ entry of $x$.  Also, ${\bf 1}_N$ denotes the
vector in $\mathbb{R}^N$, all of whose entries are equal to
one.

We shall use $V$ to denote the projection matrix in $\mathbb{R}^{N
\times N}$ given by
\begin{eqnarray}
    [V_{ij}] &:=&\begin{cases}
                    \frac{N-1}{N} & j=i\\
                    -\frac{1}{N} & j\neq i
                \end{cases}, \quad\quad i,j =1,\ldots,N,
\end{eqnarray}
and $V \mathbb{R}^N$ shall denote the image of
$\mathbb{R}^N$ under $V$.  Formally,  \[ V
\mathbb{R}^N:=\{x\in \mathbb{R}^N:~\sum_{j=1}^N x_j =
0\}.\]

\subsection{The Basic Kuramoto Model}
The basic Kuramoto model of phase-coupled oscillators under the
assumption of all-to-all coupling is given by
\begin{eqnarray}\label{eq:kuramotomodel}
    \dot{\theta_i} & = & \omega_i + \frac{k}{N}
    \sum_{j=1}^N\sin(\theta_j-\theta_i),\quad i =
    1,\ldots, N.
\end{eqnarray}
Here, $\theta_i(\cdot)\in \mathbb{R}$ ($S^1$) and
${\omega}_i \in \mathbb{R}$ respectively denote the phase
and intrinsic (or natural) frequency of oscillator~$i$, and
the constant $k\in \mathbb{R}_{+}$ is a global coupling
coefficient.

This model can be described more compactly in vector notation as
\begin{eqnarray}\label{eq:kurmodelvec}
    \dot{\theta} & = & \omega + kf(\theta)
\end{eqnarray}
where $\theta(t) := (\theta_1(t), \ldots, \theta_N(t))$,
$\omega := (\omega_1, \ldots , \omega_N)$, and the mapping
$f:\mathbb{R}^N \mapsto \mathbb{R}^N$ is given by
\begin{eqnarray}\label{eq:defoffstandard}
\nonumber f(\xi) &=& (f_1(\xi), \ldots, f_N(\xi)), \\
{f}_i(\xi) &:=& \frac{1}{N}\sum_{j=1}^{N}\sin(\xi_j-\xi_i) \quad\;\;
1 \leq i \leq N.
\end{eqnarray}
The assumption of all-to-all coupling is naturally very restrictive,
and ought to be relaxed in order for this work to be more directly
applicable to the modelling of biological systems, or most
engineering systems for that matter. Work towards this end is
underway and we hope to be able to present some results in the near
future. Meanwhile, in this paper, we shall focus exclusively on
configurations with all-to-all coupling.
%
First of all, we recall some fundamental notions in the theory of
synchronized oscillators.

\subsection{The order parameter}
\begin{figure}[t]
      \centering
      \psfrag{R}[ct][ct][1][0]{$R(\theta)$}
      \psfrag{P}[bl][bl][1][0]{$\psi(\theta)$}
      \includegraphics[width=4cm]{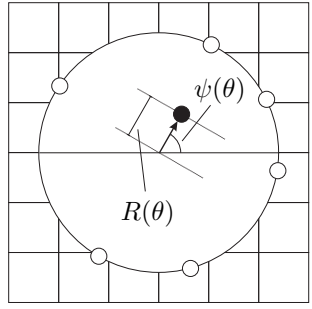}
      \caption{The order parameter $r(\theta):=R(\theta) e^{i\psi(\theta)}$ is defined as the centroid (closed circle) of the set of unit vectors (open circles) associated with the phases of the oscillators.}
      \label{fig:orderparameter}
 \end{figure}

Let $\mathbb{D}$ denote the complex unit disc $\{z\in
\mathbb{C}:~|z|\leq 1\}$. Then define $r:\mathbb{R}^N\mapsto
\mathbb{D}$, by:
\begin{eqnarray}\label{eq:orderparameter}
r(\xi) & := & \frac{1}{N}\sum_{j=1}^Ne^{i \xi_j}.
\end{eqnarray}
Let $r^{-1}(z):=\{\xi \in \mathbb{R}^N:~r(\xi)=z\}$ denote the
preimage of $r$, and note that the preimage is nonempty for all
$z\in \mathbb{D}$ provided $N\geq 2$. We introduce the notation
$\mathcal{R}_0:=r^{-1}(0)$.  Then, for $\xi \in \mathbb{R}^N$, we
may express $r(\xi)$ in polar coordinates:
\begin{equation}\label{eq:polar} r(\xi) = \begin{cases}
R(\xi)e^{i\psi(\xi)} & \xi \in
\mathbb{R}^N\backslash\mathcal{R}_0 \\
0 & \xi \in \mathcal{R}_0\end{cases}.\end{equation} Here,
$R:\mathbb{R}^N\mapsto [0,1]$ and $\psi:\mathbb{R}^N\backslash
\mathcal{R}_0 \mapsto [0,2\pi)$ are respectively defined as
\begin{eqnarray}\label{eq:absolutevalue}
    R(\xi):= \sqrt{
    \left(\frac{1}{N}\sum_{j=1}^N\sin(\xi_j)\right)^2 +
    \left(\frac{1}{N}\sum_{j=1}^N\cos(\xi_j)\right)^2},
\end{eqnarray}
and
\begin{eqnarray}\label{eq:phase}
    \psi(\xi) &:=& \arctan\left(
    \frac{\frac{1}{N}\sum_{j=1}^N\sin(\xi_j)}{\frac{1}{N}\sum_{j=1}^N\cos(\xi_j)}\right).
\end{eqnarray}
The following properties of the maps $R(\cdot)$ and $\psi(\cdot)$
follow immediately from Eqn.~(\ref{eq:orderparameter}):
\begin{equation}
    R(\xi+c{\bf 1}_N) := \left|\frac{1}{N}\sum_{j=1}^Ne^{i(\xi_j + c)}\right|
    =|e^{ic}|R(\xi) = R(\xi) \;\; \forall \xi \in \mathbb{R}^N;
\end{equation}
\begin{equation}
\psi(\xi + c{\bf 1}_N) = \psi(\xi) + c \mod 2\pi \;\; \forall \xi
\in \mathbb{R}^N\backslash \mathcal{R}_0.
\end{equation}
In the physics literature, $r(\cdot)$ is known as the \emph{order
parameter}, and is used to characterize the amount of order or
synchronization in the system (\ref{eq:kuramotomodel}). The idea is
to think of the phase $\theta_j$ of oscillator $j$ as a unit vector
$e^{i\theta_j}$ in $\mathbb{C}$; the order parameter then
corresponds to the geometric centroid of the set of vectors
$\{e^{i\theta_j}:~j=1,\ldots, N\}$, as illustrated in
Figure~\ref{fig:orderparameter}. The magnitude of the order
parameter, given by $R(\theta)$, serves as a measure of the order in
the system, in the sense that the closer the vectors are to being
perfectly aligned, the closer $R(\theta)$ is to its maximal value 1,
while vectors that are far from alignment will give rise to values
of $R(\theta)$ significantly smaller than 1.

It follows from (\ref{eq:polar}) that for $\xi \in
\mathbb{R}^N\backslash \mathcal{R}_0$,
\begin{eqnarray}\label{eq:Rtheta1} R(\xi) &=& e^{-i
\psi(\xi)}r(\xi)\\ &=& \frac{1}{N} \sum_{j=1}^{N} e^{i\left(\xi_j -
\psi(\xi)\right)}.
\end{eqnarray}
Equating real and imaginary parts in (\ref{eq:Rtheta1}), we
immediately see that for $\xi \in \mathbb{R}^N\backslash
\mathcal{R}_0$:
\begin{equation}\label{eq:usefulidentity} R(\xi) =
\frac{1}{N}\sum_{j=1}^N\cos(\psi(\xi)-\xi_j);\end{equation}
\begin{equation}\label{eq:usefulidentity2} \sum_{j=1}^N\sin(\psi(\xi)-\xi_j) =
0.\end{equation} Both of these identities shall prove useful
throughout the paper.

Before proceeding, note that the function $f:\mathbb{R}^N \mapsto
\mathbb{R}^N$ given by (\ref{eq:defoffstandard}) can be written in
terms of the functions $R(\cdot)$ and $\psi(\cdot)$ as
\begin{equation}\label{eq:defoff}
    f_i(\xi) := \begin{cases}R(\xi)\sin(\psi(\xi)-\xi_i) & \xi\in
    \mathbb{R}^N\backslash \mathcal{R}_0 \\
    0 & \xi\in \mathcal{R}_0\end{cases},
\end{equation}
for $1 \leq i \leq N$.

\subsection{Fixed points and global phase-locking}
Let $\langle\omega\rangle$ denote the sample mean of the natural
frequencies,
$\langle\omega\rangle:=\frac{1}{N}\sum_{j=1}^N\omega_j$. Similarly,
let
 $\langle \theta(t) \rangle$ denote the mean phase of a solution of (\ref{eq:kuramotomodel}) at time~$t$.  In general, $\langle \omega \rangle$ and $\langle \theta(t)
 \rangle$ will be non-zero.  However, we shall now show that for the
 study of phase-locked solutions of (\ref{eq:kuramotomodel}), we may assume without loss of
 generality that $\langle \omega \rangle = 0$, $\langle \theta(t)
 \rangle = 0$ for $t \geq t_0$.  This helps to simplify the analysis
 of phase-locked solutions of (\ref{eq:kuramotomodel}), as it allows
 us to transform the problem into a question of fixed point
 existence for a lower-dimensional system.

Consider the new coordinates
 \begin{eqnarray}
    x_i(t) &:=& \theta_i(t) - \langle \theta(t) \rangle, \quad
    i=1,\ldots,N.
 \end{eqnarray}
Then $x(t):=V\theta(t)$. Similarly, define
$\Omega:=V\omega$. In the new coordinates, the system
dynamics are given by:
\begin{equation}\label{eq:kuramotomodelgrounded} \dot{x} = \Omega +
kf(x), \quad\; x(t) \in V\mathbb{R}^N,
\end{equation}where $f(\cdot)$ is defined in
(\ref{eq:defoffstandard}).

The key point here is that as $\langle \Omega \rangle = 0$,
$V \mathbb{R}^N$ is invariant under
(\ref{eq:kuramotomodelgrounded}). To avoid confusion, we
shall use $x(t)$ to denote solutions to the system
(\ref{eq:kuramotomodelgrounded}) on $V \mathbb{R}^N$, while
$\theta(t)$ shall be used to denote solutions to the
original Kuramoto system (\ref{eq:kuramotomodel}) on
$\mathbb{R}^N$. Our main concern for the remainder of the
paper is to find conditions on $k$ and $\Omega$ under which
the system~(\ref{eq:kuramotomodelgrounded}) has one or more
fixed points in the sense of the following definition.
\begin{defn}[fixed point]\label{def:fixedpoint}Given $\omega\in \mathbb{R}^N$, let $\Omega:=V\omega$.  We say that $x \in V\mathbb{R}^N$
is a \emph{fixed point} (of the system (\ref{eq:kuramotomodelgrounded})) if
\begin{equation}\label{eq:fixedpointequation}
    kf(x) = -\Omega.
\end{equation}
\end{defn}

There is a natural correspondence between fixed points of
(\ref{eq:kuramotomodelgrounded}) and phase-locked solutions of
(\ref{eq:kuramotomodel}).  In fact, for every fixed point $x^*\in
V\mathbb{R}^N$ there is a $1$-dimensional manifold
$\mathcal{M}:=\{\theta \in \mathbb{R}^N:~\theta = x^*+\langle \omega
\rangle t,~t\in \mathbb{R}\}$ that is invariant under the original
system dynamics~(\ref{eq:kuramotomodel}). More precisely, let
${x}^*$ be a {fixed point} and let $\theta^0\in \mathbb{R}^N$ be
such that $V\theta^0 = x^*$. Then the solution $\theta(t)$ of the
system~(\ref{eq:kuramotomodel}) with initial condition
$\theta(t_0)=\theta^0$ satisfies
\begin{eqnarray}
    \theta_i(t) - \theta_j(t) = \theta^0_i - \theta^0_j
\end{eqnarray}
for all $t\geq t_0$ and all $(i,j)$. In other words, a fixed point
in the sense of Definition~\ref{def:fixedpoint} corresponds to a
situation in which each oscillator is phase-locked to every other
and moves at constant speed $\dot{\theta}_i = \langle \omega
\rangle$. We shall refer to this phenomenon as global
phase-locking. In the literature, it is also known as full (or
complete) synchronization. See also~\cite{Jad04}.
\subsection{Critical coupling}
We next define the notion of {\it critical coupling }, which is
central to the work of the rest of the paper.  Essentially, the
critical coupling is the smallest $k$ for which the
system~(\ref{eq:kuramotomodelgrounded}) has at least one fixed
point.  Formally, we have the following definition.
\begin{defn}\label{def:criticalcoupling}Given $\omega \in \mathbb{R}^N$, let $\Omega := V\omega$. We define the \emph{critical
coupling}, $k_\text{c}$, as follows:
\begin{eqnarray}
    k_\text{c} &:=& \inf_k \left\{k\in \mathbb{R}_{+}: \exists x\in
V\mathbb{R}^N \text{s.t.}~kf(x) = -\Omega\right\}.
\end{eqnarray}
\end{defn}
Note that this definition of the critical coupling, which is
equivalent to that of $K_L$ in \cite{Jad04}, does not coincide
with the traditional notion used in the physics literature.
Indeed, the traditional notion of critical coupling is defined in
terms of the lowest value of $k$ for which there exists at least
one solution $x(t)$, $t\geq t_0$, and a constant $c\in (0,1]$ such
that $R(x(t))=c$ for all $t\geq t_0$ (so called stationary or
steady solutions \cite{Stro00}). Note that these solutions are not
necessarily fixed points, although, in finite dimensions, the
probability of finding stationary solutions that are not fixed
points is vanishingly small.  In his original analysis, Kuramoto
showed that in the limiting case when $N$ tends to infinity,
stationary solutions always exist for large enough $k$, provided
the distribution of natural frequencies is symmetric. Our
definition, although more restrictive in a sense, does not impose
any restriction on the shape of the distribution of natural
frequencies other than that it should have compact support.  In
fact, it follows from the result of Lemma \ref{lem:lowerbound}
below that, if the distribution of natural frequencies does not
have compact support, then the critical coupling will exceed any
finite number with probability tending to $1$ as $N$ tends to
infinity. In this paper we shall therefore focus on distributions
with compact support.

\section{Existence of fixed points under strong coupling}
In this section we shall show that, provided the distribution of
intrinsic frequencies has compact support, the critical coupling
given in Definition~\ref{def:criticalcoupling} is always finite.
Following on from this, in the next section, we shall derive a
number of lower bounds for the value of the critical coupling.
There are two steps in the derivation given here: first, we
characterize the fixed points of the homogeneous system
\begin{equation}\label{eq:homogenous}\dot{x} = kf(x).\end{equation} Then, in the second step, we use
a perturbation argument to show that for every fixed point of the
homogeneous system, we can find an open set containing it, such
that, under strong enough coupling, the original
system~(\ref{eq:kuramotomodelgrounded}) has a unique fixed point on
this set. As a first step, the following lemma characterises the
fixed points of the homogeneous system.

\begin{lemma}\label{lem:fixedpointnecsuf}  Let $f(\cdot)$ be given by
(\ref{eq:defoffstandard}) and $\xi \in \mathbb{R}^N$. We
have that $f(\xi)=0$ if and only if one or both of the
following conditions are satisfied
\begin{quote}
\begin{itemize}
\item[$(a)$] $R(\xi) = 0$; \item[$(b)$]
$\sin(\xi_i-\xi_j)=0$ for all $(i,j)$.
\end{itemize}
\end{quote}
\end{lemma}

\begin{proof}Sufficiency of conditions $(a)$ and $(b)$ follows
from (\ref{eq:defoff}) and (\ref{eq:defoffstandard})
respectively. To prove necessity, suppose $f(\xi)=0$ and
 $R(\xi)\neq 0$ (if $R(\xi)=0$ we are done). It follows that
$\sin(\psi(\xi)-\xi_i)=0$ for all $i$. This implies that
there exist integers $k_i\in \mathbb{Z}$, $i=1,2,\ldots N$
such that $\psi(\xi)-\xi_i = k_i\pi $ for all $i$, and we
have that $\xi_i-\xi_j = (k_j-k_i)\pi$. We conclude that
$\sin(\xi_i-\xi_j)=0$ for all $(i,j)$.
\end{proof}

%
%

\begin{remark}It is not hard to see that conditions $(a)$
and $(b)$ in Lemma~\ref{lem:fixedpointnecsuf} are mutually
exclusive if and only if the dimension $N$ is odd. We shall
prove necessity. Suppose conditions $(a)$ and $(b)$ both
hold and suppose furthermore that $N$ is odd. Then for all
$(i,j)$ we have that either $\cos(\xi_i-\xi_j)=1$ or
$\cos(\xi_i-\xi_j)=-1$. We write
$R^2(\xi)=\frac{1}{N^2}\sum_{i,j}\cos(\xi_i-\xi_j)=\frac{1}{N^2}\left(N+2\sum_{i,j>i}\cos(\xi_i-\xi_j)\right)$.
Since $R(\xi)=0$ by assumption, it follows that
$2\sum_{i,j>i}\cos(\xi_i-\xi_j)=-N$. The left hand side
evaluates to an even integer. By assumption, the number on
the right hand side is odd. We arrive at a contradiction
and conclude that if $N$ is odd, conditions $(a)$ and $(b)$
cannot both hold.
\end{remark}

Next we shall prove that the fixed points of our $N$-dimensional
system~(\ref{eq:kuramotomodelgrounded}) can be found by solving a
system of $N-1$ equations in as many variables. We have the
following result:

\begin{lemma}\label{lem:reduce}
Let $p\in \{1,\ldots,N\}$ and let $x^{*}\in V\mathbb{R}^N$. Then
$x^{*}$ is a fixed point of~(\ref{eq:kuramotomodelgrounded}) if and
only if $kf_i(x^{*}) = -\Omega_i$ for $i\neq p$.
\end{lemma}

\begin{proof}
The proof of necessity is trivial. To prove sufficiency,
recall that \begin{equation}\label{eq:usid}
\sum_{j=1}^N\left(\Omega_j+kf_j(x)\right) = 0\quad
\text{for all}~x\in \mathbb{R}^N. \end{equation} Now
suppose $kf_i(x)= -\Omega_i$ for all $i\neq p$. Then it
follows from~(\ref{eq:usid}) that $\Omega_p+kf_p(x) = 0$.
In other words, it follows that $kf_i(x)=-\Omega_i$ for
\emph{all} $i$. We conclude that $x$ is a fixed point.
\end{proof}

Let $x^*\in V\mathbb{R}^N$ be a fixed point of the homogeneous
system~(\ref{eq:homogenous}) such that $R(x^*)\neq 0$. We shall now
show that \emph{locally}, in a neighborhood of $x^*$, the system of
equations
\begin{equation}\label{eq:systemofequations}
\left\{
\begin{array}{ccc}
    -{\Omega_1} & = & kf_1(x_1,\ldots,x_{N-1},-\sum_{j=1}^{N-1} x_j)
     \\
    \vdots &  & \vdots \\
    -{\Omega_{N-1}} & = &
    kf_{N-1}(x_1,\ldots,x_{N-1},-\sum_{j=1}^{N-1}x_j)
\end{array}\right.
\end{equation}
has a unique solution, provided $k$ is large enough. It
follows directly from Lemma~\ref{lem:reduce} that every
solution of~(\ref{eq:systemofequations}) defines a fixed
point and, conversely, that every fixed point
satisfies~(\ref{eq:systemofequations}). We proceed as
follows. Let $x^*\in V\mathbb{R}^N$. We define the Jacobian
$J(x^*)$, as follows:
\begin{eqnarray}\label{eq:jacobian}
    [J_{ij}(x^*)] & := & \left.\frac{\partial
    f_i(x_1,\ldots,-\sum_{j=1}^{N-1}
x_j)}{\partial
    x_j}\right|_{x=x^*},
\end{eqnarray}
where $i,j = 1,\ldots, N-1$. We have the following result:

\begin{lemma}\label{lem:nonsingular} Let $f(\cdot)$ be given by
(\ref{eq:defoffstandard}) and suppose that $x^* \in V\mathbb{R}^N$
satisfies $f(x)=0$ and $R(x)\neq 0$.  Then $\det(J(x^*))\neq 0$.
\end{lemma}

\begin{proof}
Let $x^*$ be a fixed point of the homogeneous system and
suppose $R(x^*)\neq 0$. Then by
Lemma~\ref{lem:fixedpointnecsuf}, we have that
$\sin(x^*_j-x^*_i) = 0$ for all $(i,j)$, and it follows
that
\begin{eqnarray}
\label{eq:cosdiff}\nonumber \cos(x^*_j-x^*_i) &=&
\cos((x^*_j-x^*_s)-(x^*_i-x^*_s)) \\ &=&
\cos(x^*_j-x^*_s)\cos(x^*_i-x^*_s)
\end{eqnarray}
for all $s$ and all $(i,j)$. The claim is that $J(x^*)$ is
nonsingular. To prove this, we proceed as follows. From the
definition, we have that
\begin{eqnarray}\label{eq:Jfirst}
J_{ij}(x^*) & = & \begin{cases} -\sum_{m = 1, m\neq
i}^{N-1}\cos(x^*_m-x^*_i) - 2\cos(x^*_N - x^*_i) & i = j \\
\cos(x^*_j - x^*_i) - \cos(x^*_N - x^*_i) & i \neq j
\end{cases}.
\end{eqnarray}
Using the aforementioned identity, setting $s=N$, we
rewrite~(\ref{eq:Jfirst}), as follows,
\begin{eqnarray}
J_{ij}(x^*) & = & \begin{cases} -\left(\sum_{m=1, m\neq
i}^{N-1}\cos(x^*_m-x^*_N) + 2\right)\cos(x^*_i - x^*_N) & i = j\\
\left(\cos(x^*_j - x^*_N) -  1\right)\cos(x^*_i-x^*_N) & i
\neq j
\end{cases}.
\end{eqnarray}
Inspection shows that the rank of $J(x^*)$ is invariant
under permutations of the components of $x^*$. Hence we can
assume, without loss of generality, that there exists
$\rho\in \{0,\ldots,N-1\}$, such that
\begin{eqnarray}\label{eq:satisfythis}
\cos(x^*_j-x^*_N) & = &
\begin{cases}
-1, & \quad 1\leq i \leq \rho \\
+1, & \quad \rho+1 \leq i \leq N
\end{cases}
\end{eqnarray}
Under this assumption $J(x^*)$ takes the form
\begin{eqnarray}
J(x^*) & = &\begin{pmatrix}A & 0 \\ C & B\end{pmatrix}
\end{eqnarray}
where $A$ and $B$ are square matrices of dimension $\rho\times
\rho$ and $(N-1-\rho) \times (N-1-\rho)$ respectively. It follows
that $J(x^*)$ is nonsingular if and only if $A$ and $B$ are
nonsingular. Inspection shows that $A=(N-2\rho)I + 2 {\bf 1}{\bf
1}^T$ and $B=(2\rho - N)I$. It follows that $A$ or $B$ is singular
if and only if $N=2\rho$. In case $N$ is odd, this condition is
never satisfied. In case $N$ is even this condition, combined with
(\ref{eq:cosdiff}) and the fact that
$R^2(x)=\frac{1}{N^2}\sum_{i,j}\cos(x_i-x_j)$, implies that
$R(x^*)=0$, which contradicts our starting assumption. We conclude
that, under the hypotheses of the lemma, $J$ is nonsingular. This
concludes the proof.
\end{proof}

Let $\Pi: \mathbb{R}^{N-1} \mapsto V\mathbb{R}^N$, be given
as
\begin{eqnarray}
(\Pi(y))_i & := & \begin{cases}y_i & \text{for}~
i=1,2,\ldots N-1;
\\ -\sum_{j=1}^{N-1} y_j & \text{for}~i=N.\end{cases}
\end{eqnarray}
and note that $\Pi$ has an inverse $\Pi^{-1}$ that is
defined everywhere in $V\mathbb{R}^N$. We are now ready to
state the main result:
\begin{theorem}\label{thm:solution}
Let $f(\cdot)$ be given by (\ref{eq:defoffstandard}) and $x^*\in V
\mathbb{R}^N$ be such that $f(x^*)=0$ and $R(x^*)\neq 0$. Also, let
$\Omega\in V\mathbb{R}^N$. Then there exists $K\in \mathbb{R}$, and
an open set $U\in \mathbb{R}^{N-1}$ such that (a) $\Pi^{-1}(x^*)$ is
an interior point of $U$; and (b) for all $k>K$, the system of
equations~(\ref{eq:systemofequations}) has a unique solution on $U$.
\end{theorem}

\begin{proof}
Define $y^*:=\Pi^{-1}(x^*)$, and let
$g:\mathbb{R}^{N-1}\mapsto \mathbb{R}^{N-1}$ be given as
$g_i(y) := f_i(y_1,\ldots, y_{N-1}, -\sum_{j=1}^{N-1}y_j)$,
$i=1,\ldots,N-1$. Note that $g(y^*)=0$. Also, by
Lemma~\ref{lem:nonsingular}, we have that
$\det(\frac{\partial g}{\partial y}(y^*))|\neq 0$. Under
these conditions, the Inverse Function Theorem says that
there exists an open set $U\subset \mathbb{R}^{N-1}$
containing $y^*$ such that $\left.g\right|_U:U\mapsto g(U)$
is a diffeomorphism. By continuity (and bijectivity) of
$g^{-1}$ there exists $\delta>0$ such that for all $z\in
\mathbb{R}^{N-1}$ satisfying $\|z\|< \delta$, the equation
$g(y)=z$ has a unique solution on $U$. Now let $z$ be given
as $z_i:=-\Omega_i/k$. Since, by assumption, $\max_i
|\Omega_i|<\infty$, it follows that, provided $k$ is large
enough, the system of equations $\{kg_i(y)=-{\Omega}:
i=1,\ldots, N-1\}$ has a unique solution on $U$. This
concludes the proof.
\end{proof}

As alluded to earlier, there is a unique correspondence
between solutions of~(\ref{eq:systemofequations}) and the
fixed points of the
system~(\ref{eq:kuramotomodelgrounded}). Indeed, by
Lemma~\ref{lem:reduce} we have that if $y$ is a solution
of~(\ref{eq:systemofequations}), then $\Pi^{-1}(y)$ is a
fixed point, and conversely, if $y$ is a fixed point, then
$\Pi(y)$ is a solution of~(\ref{eq:systemofequations}).
Thus, an immediate consequence of
Theorem~\ref{thm:solution} is that for large enough $k$,
the system~(\ref{eq:systemofequations}) will have at least
one fixed point. In other words, Theorem~\ref{thm:solution}
tells us that the critical coupling, $k_\text{c}$, is
always finite.

Note furthermore that the proof of
Theorem~\ref{thm:solution} does not require detailed
knowledge of the coupling function $g$ and that, as such,
its applicability is not restricted to networks with
all-to-all coupling. To illustrate this, consider the case
of a $4$-cycle, where $g$ is given as
\begin{eqnarray}
    g_1(y) & = & \frac{1}{4}\sin(y_2-y_1) +
    \frac{1}{4}\sin(-2y_1-y_2-y_3)\nonumber \\
    g_2(y) & = & \frac{1}{4}\sin(y_3-y_2) +
    \frac{1}{4}\sin(y_1-y_2)\nonumber \\
    g_3(y) & = & \frac{1}{4}\sin(-2y_3-y_1-y_2) +
    \frac{1}{4}\sin(y_2-y_3).
\end{eqnarray}
We have that $g(0)=0$ and $\det(\frac{\partial g}{\partial
y}(0))=-\frac{1}{4}\neq 0$. This implies that for $k$ large
enough, the system of equations $\{k g_i(y) = -\Omega_i:
i=1,2,3\}$ has a unique solution on some open set
containing the origin.

Lastly, note that continuity of $g^{-1}$ implies that the
fixed points of the original
system~(\ref{eq:kuramotomodelgrounded}) will converge to
the fixed points of the homogeneous
system~(\ref{eq:homogenous}) as $k$ tends to infinity.

\section{Lower bounds on the critical coupling\label{sec:lowerbounds}}
In the previous section we showed that the critical
coupling is finite, provided the oscillator's intrinsic
frequencies are finite. In the present section we shall
investigate in more detail the relation between the
distribution of intrinsic frequencies and the critical
coupling. In particular, we shall derive various lower
bounds and discuss some of these bound's implications for
the system's dynamic behaviour.

First, let us observe that $k_\text{c}$ (Definition
\ref{def:criticalcoupling}) is lower bounded by the
$l^{\infty}$ norm of $\Omega$:
\begin{equation}\label{eq:lowerboundsimple}
    k_\text{c} \; \geq \; \|\Omega\|_{\infty} \; := \;\max_i |\omega_i - \langle
    \omega
    \rangle|.
\end{equation}
This follows trivially from inspection of
Eqn.~(\ref{eq:kuramotomodelgrounded}). In order to derive
another lower bound, we shall need the following result:
\begin{lemma}\label{lem:lowerbound}
    Let $f(\cdot)$ be given by~(\ref{eq:defoffstandard}). Then:
    \begin{enumerate}
    \item For all $x \in
    \mathbb{R}^N$,
    \begin{eqnarray}\label{eq:inequality}
        \|f(x)\|_{2} & \leq &
        \sqrt{NR^2(x)\left(1-R^2(x)\right)};
    \end{eqnarray}
    \item If $N$ is even, then for every $c \in [0,1]$ there exists
    $x\in V\mathbb{R}^N$ such that $R(x)=c$ and
    $\|f(x)\|_{2}=\sqrt{NR^2(x)\left(1-R^2(x)\right)}$;
    \item If
    $N$ is odd, then inequality (\ref{eq:inequality}) is strict for all $x\in \mathbb{R}^N$
    such that $0<R(x)<1$.
    \end{enumerate}

    \end{lemma}
\begin{proof}
    Part 1. Observe that inequality~(\ref{eq:inequality}) is trivially satisfied when $x\in\mathcal{R}_0$. Suppose therefore that $x\in \mathbb{R}^N\backslash
    \mathcal{R}_0$. Then by definition
    \begin{eqnarray}\label{eq:step1}
        \|f(x)\|^2_{2} &:=&
        \sum_{j=1}^N(f_j(x))^2\nonumber \\
        & = &  R^2(x) \sum_{j=1}^N
        \sin^2(\psi(x)-x_j),
    \end{eqnarray}
where $\psi(x)$ and $R(x)$ are the phase and magnitude of
the order parameter, previously defined in (\ref{eq:phase})
and (\ref{eq:absolutevalue}) respectively. Introducing the
shorthand notation $z_i(x):=\cos(\psi(x)-x_i)$, and using
(\ref{eq:usefulidentity}) we now rewrite~(\ref{eq:step1}),
as follows:
\begin{eqnarray}
    \|f(x)\|^2_{2} & = & \left(\frac{1}{N}\sum_{j=1}^N z_j(x)\right)^2
    \sum_{j=1}^N
        \left(1-z_j(x)^2\right).
\end{eqnarray}
To derive the desired inequality we pick a $c\in [0,1]$ and
maximize $\|f(x)\|_2$ over the set $\{x\in
\mathbb{R}^N:R(x)=c\}$. We shall not solve this
optimization problem directly, but take an indirect route
by considering another, easier optimization problem, whose
solution will then give us an upper bound on the solution
to the first problem. Then we shall show that, under
certain conditions, the two solutions coincide.

To this end, let $c\in (0,1]$ and consider the constrained
optimization problem
\begin{center}\parbox{1.5cm}{OPT 1:}\framebox[7.3cm][c]{\begin{minipage}[c]{7.1cm}\vspace{-0.25cm}
\begin{equation}
\begin{array}{lll}\text{maximize} &
 \sum_{j=1}^N\left(1-z_j(x)^2\right)&\\
  \text{subject to} & \frac{1}{N}\sum_{j=1}^N z_j(x) =
c, & {x \in \mathbb{R}^N\backslash
\mathcal{R}_0}\end{array}\nonumber\end{equation}
\end{minipage}}
\end{center}
Note that the constraint is feasible for all values of $c$
in the specified interval. We shall denote the solution to
OPT 1 as $s_1(c)$. Next consider a second optimization
problem,
\begin{center}\parbox{1.5cm}{OPT 2:}\framebox[7.3cm][c]{\begin{minipage}[c]{7.1cm}\vspace{-0.25cm}
\begin{equation}
\begin{array}{lll}\text{maximize} & \sum_{j=1}^N\left(1-y_j^2\right) & \\
  \text{subject to} & \frac{1}{N}\sum_{j=1}^N y_j =
c, & {y \in \mathbb{R}^N}.
 \end{array}\nonumber\end{equation}\end{minipage}
} \end{center} and let the solution to this problem be
denoted as $s_2(c)$. We then have that $s_2(c)\geq s_1(c)$
for all $c\in (0,1]$. In other words, the solution to OPT 1
is upper bounded by the solution to OPT 2. The solution to
OPT 2 can be found by means of standard Lagrange multiplier
techniques. The optimum $s_2(c)=N\left(1-c^2\right)$, is
attained when $y_i=c$ for all $i$. We conclude that
\begin{equation}
\max_{\{x\in \mathbb{R}^N:R(x)=c\}}\|f(x)\|^2_2 \leq
Nc^2\left(1-c^2\right),
\end{equation}
and hence,
\begin{equation}
   \|f(x)\|_2 \leq \sqrt{N}R(x)\sqrt{1-R^2(x)}.
\end{equation}
for all $x\in \mathbb{R}^N$.

Part 2. To prove the second part of the theorem, let $c\in
(0,1]$ and note that $s_1(c) = s_2(c)$ if and only there
exists $x \in \mathbb{R}^N\backslash \mathcal{R}_0$ such
that
\begin{equation}\label{eq:conditionforoptimality}
 \cos(\psi(x)-x_i) = c
\end{equation}
for all $i$. Suppose $N$ is even and let $x$ be given as
\begin{eqnarray}
    x_i & := & \begin{cases}\arccos(c) & i = 1,\ldots,
    \frac{N}{2} \\
    -\arccos(c) & i = \frac{N}{2},\ldots,N.
                \end{cases}
\end{eqnarray}
Then $\sum_{j=1}^Nx_j=0$, and, by definition, $x\in
V\mathbb{R}^N$. Moreover, $\psi(x)=0$, and
$\cos(\psi(x)-x_i) = c$ for all $i$. This completes the
second part.

Part 3. To prove the third part, let $N$ be odd and suppose
there exists $x\in \mathbb{R}^N$ such that
Condition~(\ref{eq:conditionforoptimality}) is satisfied.
Then it follows from the identity
$\sin^2(\psi(x)-x_i)+\cos^2(\psi(x)-x_i)=1$ that there must
exist $a\in \{-1,1\}^N$ such that $\sin(\psi(x)-x_i) =
a_i\sqrt{1-c^2}$ for all $i$. By
Identity~(\ref{eq:usefulidentity2}), we have that $\sum_j
\sin(\psi(x)-x_j)=0$, which, assuming $c\neq 1$, implies
that $\sum_{j=1}^N a_j=0$. But this cannot be true unless
$N$ is even. Thus we arrive at a contradiction and we
conclude that if $N$ is odd then $s_2(c)>s_1(c)$ for all
$c$ such that $0<c<1$. This concludes the proof.
\end{proof}

\begin{figure}[thpb]
\centering

\psfrag{G}[ct][ct][0.8][0]{$0$}
\psfrag{H}[ct][ct][0.8][0]{$0.2$}
\psfrag{I}[ct][ct][0.8][0]{$0.4$}
\psfrag{J}[ct][ct][0.8][0]{$0.6$}
\psfrag{K}[ct][ct][0.8][0]{$0.8$}
\psfrag{L}[ct][ct][0.8][0]{$1$}

\psfrag{A}[cr][cr][0.8][0]{$0$}
\psfrag{B}[cr][cr][0.8][0]{$0.1$}
\psfrag{C}[cr][cr][0.8][0]{$0.2$}
\psfrag{D}[cr][cr][0.8][0]{$0.3$}
\psfrag{E}[cr][cr][0.8][0]{$0.4$}
\psfrag{F}[cr][cr][0.8][0]{$0.5$}

\psfrag{R}[cr][cr][0.8][0]{$R(x)\rightarrow$}
\psfrag{O}[ct][ct][0.8][90]{$\frac{1}{\sqrt{N}}\|f(x)\|_{2}\rightarrow$}

    \begin{tabular}{cc}

      \includegraphics[width=4cm]{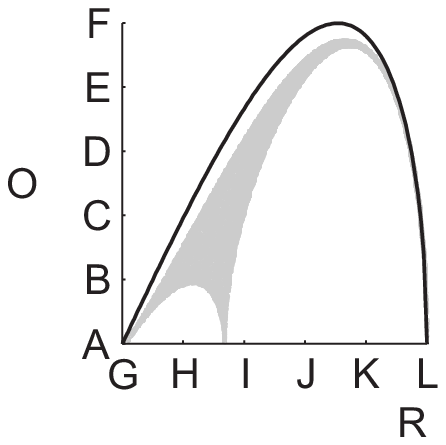}
      &
      \includegraphics[width=4cm]{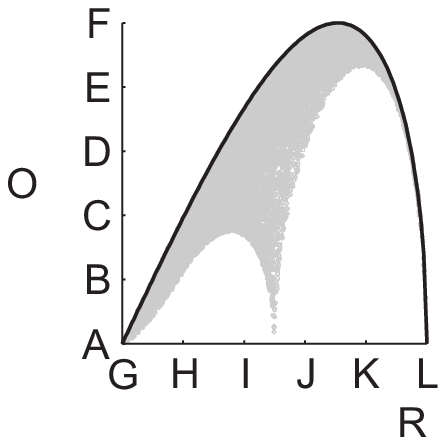}\\
      (a) $N=3$ & (b) $N=4$

    \end{tabular}

      \caption{Scatter plot of
      $\frac{1}{\sqrt{N}}\|f(x)\|_2$ for $N=3$ (left panel) and $N=4$ (right panel). The phases $x$ were drawn from a uniform distribution. The solid black line in both panels is the upper bound
      $R(x)\sqrt{1-R^2(x)}$.}
      \label{fig:scatterplots}
 \end{figure}

Figure~\ref{fig:scatterplots} illustrates the result of
Lemma~\ref{lem:lowerbound}. When $N=4$ (even), the lower
bound is attained at every value of $R(x)$, which shows
that the given bound is the tightest possible. However, as
illustrated in the left panel, when $N=3$, the bound is
never attained except on the set $\{x\in
\mathbb{R}^N:R(x)\in \{0,1\}\}$. It can be shown however
that in the limit of large $N$ the given bound is
arbitrarily tight, even for odd $N$, in the sense that for
every $c\in [0,1]$,
\[
\min_{\{x\in
\mathbb{R}^{2m+1}:R(x)=c\}}\frac{1}{\sqrt{2m+1}}\left|\|f(x)\|_2
- \sqrt{(2m+1)c^2\left(1-c^2\right)} \right|\nonumber
\]
tends to zero as $m$ tends to infinity.

Lemma~\ref{lem:lowerbound} has some interesting
implications. For instance, it can provide insight into the
rate at which solutions of a homogeneous system of Kuramoto
oscillators ((\ref{eq:kuramotomodel}) with $\omega_i = 0$
for $1 \leq i \leq N$) on $\mathbb{R}^N$ converge to fixed
points, . To see this, consider the homogeneous system
\begin{equation}\label{eq:homsystem}
\left\{\begin{array}{ccc}\dot{\theta}(t) & = & kf(\theta(t))\\
\theta(t_0) & = &\theta_0\end{array}\right.,
\end{equation}
where $\theta_0\in \mathbb{R}^N$. We shall compute the time
derivative of the magnitude squared of the order parameter,
$L(\cdot):=R^2(\cdot)$, and show that this derivative is: $(i)$
nonnegative along solutions of~(\ref{eq:homsystem}); ($ii$) bounded
from above by a certain function $D(t)$ for every $t$. We proceed as
follows~\cite{Jad04}. By definition,
\[
\frac{\text{d}L(\theta(t))}{\text{d}t} := \frac{L(\theta)}{\partial
\theta}\dot{\theta}(t) = \frac{L(\theta)}{\partial
\theta}kf(\theta(t)).
\]
Using the identity
\[
\frac{\partial L(\theta)}{\partial \theta} =
\frac{2}{N}\left[f(\theta)\right]^T,
\]
it follows that
\begin{equation}\label{eq:Rderivative}
\frac{\text{d}L(\theta(t))}{\text{d}t} = \frac{2k}{N}
\|f(\theta(t))\|^2_2,
\end{equation}
which shows that the time-derivative is positive
everywhere, except at the equilibria, where it is zero. It
follows that the magnitude of the order parameter is a
nondecreasing function of time. Based on the observation
that the time-derivative of $L$ is positive almost
everywhere (the set of equilibria having measure zero), we
formulate the following conjecture, \cite{Cho05,Jad04}:

\begin{conjecture}\label{con:globalstable}
For almost all initial conditions $\theta_0$, the solution
$\theta(t)$ to the homogeneous system~(\ref{eq:homsystem}) has the
property that $\lim_{t\rightarrow \infty} R(\theta(t))=1$.
\end{conjecture}

In agreement with Conjecture~\ref{con:globalstable}, one can prove
that, for the homogeneous system, the global phase-locking manifold
$\mathcal{M}:=\{\theta \in \mathbb{R}^N: \theta_i =
\theta_j~\text{for all}~i,j\}$ is (locally) asymptotically stable.
However, the existence of other invariant manifolds, not contained
in $\mathcal{M}$, implies that $\mathcal M$ is not globally
asymptotically stable. We conjecture that $\mathcal{M}$ is `almost
globally asymptotically stable', in the sense that its region of
attraction is the entire space minus a set of measure zero.

For our next result, we shall need the concept of a
dominating function, which is defined as follows:
\begin{defn}Let $h,g:\mathbb{R}\mapsto \mathbb{R}$ and let
$\mathcal{I} \subset \mathbb{R}$ be some interval. We say that $h$
dominates $g$ on $\mathcal{I}$ if $h(t)\geq g(t)$ for all $t\in
\mathcal{I}$. In that case we call $h$ a \emph{dominating function}
for $g$ on $\mathcal{I}$.
\end{defn}

Our next result states that $L(\theta(t))$ is dominated by
a certain scalar function $D(t)$ that depends only on
$\theta_0$. In order to prove this result, we need the
following two lemmas.

\begin{lemma}\label{lem:equalalltimes}Let $\theta(\cdot)$ be a solution of the homogeneous system~(\ref{eq:homsystem}) and suppose $\dot{L}(\theta(t')) =
2k L(\theta(t'))\left(1-L(\theta(t'))\right)$ for some
$t'\in \mathbb{R}$. Then
\[
\dot{L}(\theta(t)) =
2kL(\theta(t))\left(1-L(\theta(t))\right)\quad \text{for
all $t\geq t'$}.
\]
\end{lemma}
\begin{proof}
Recall that $\dot{L}(\theta(t))=\frac{2k}{N}\|f(\theta(t))\|^2_2$.
It follows from the proof of Lemma~\ref{lem:lowerbound} that
$\|f(\theta(t'))\|^2_2 = N
L(\theta(t'))\left(1-L(\theta(t'))\right)$ for some $t'\in
\mathbb{R}$ if and only if one or two of the following conditions
hold: (a) $L(\theta(t'))=0$; (b) $N$ is even and there exists a
permutation $\hat{\theta}(t')$ of $\theta(t')$ such that
\begin{eqnarray}\label{eq:conditiona}
\cos(\hat{\theta}_i(t') - \hat{\theta}_1(t')) & = &
1,\quad\quad i=1,2,\ldots,\frac{N}{2} \\
\label{eq:conditionb} \cos(\hat{\theta}_i(t') -
\hat{\theta}_N(t')) & = & 1,\quad\quad i=\frac{N}{2} +
1,\ldots, N
\end{eqnarray}
If $L(\theta(t'))=0$, we have that $L(\theta(t))=0$ for all
$t\geq t'$ and the result follows trivially. Now suppose
Conditions (\ref{eq:conditiona}) and (\ref{eq:conditionb})
hold. It follows that $\dot{\hat{\theta}}_i(t') =
\dot{\hat{\theta}}_j(t')$ for $i,j \leq \frac{N}{2}$ and
$i,j
> \frac{N}{2}$, and hence
\[
\frac{\text{d}}{\text{d}\tau}\left.\left(\cos(\hat{\theta}_i(\tau)-\hat{\theta}_j(\tau))\right)\right|_{\tau=t'}
= 0, \quad i,j \leq \frac{N}{2}\; ; \;i,j > \frac{N}{2}
\] In other words, if $\hat{\theta}(\cdot)$
satisfies Conditions (\ref{eq:conditiona}) and
(\ref{eq:conditionb}) for some $t'$, it
satisfies~(\ref{eq:conditiona}) and (\ref{eq:conditionb})
 for all $t\geq t'$. This concludes the proof.
\end{proof}

\begin{lemma}\label{lem:dominator} Let $h,g:\mathbb{R} \mapsto \mathbb{R}$ be
such that for every $x_0\in \mathbb{R}$, the systems
\[\left\{\begin{array}{ccc}\dot{x} & = & h(x) \\ x(0) & = &
x_0\end{array}\right., \quad \quad
\left\{\begin{array}{ccc}\dot{x} & = & g(x) \\ x(0) & = &
x_0\end{array}\right.
\]
have unique solutions in $C^1[0,\infty)$. Let these solution be
denoted $x_h(t;x_0)$ and $x_g(t;x_0)$, respectively. Suppose there
exists $a,b\in \mathbb{R}$, $a\neq b$, such that \[h(x) > g(x)\geq 0
\] for all $x\in (a,b)\subset \mathbb{R}$. Then for every $x_0\in
(a,b)$, we have that
\begin{eqnarray}
x_h(t;x_0) \geq x_g(t;x_0)
\end{eqnarray}
for all $t\in I$, where $I\subset \mathbb{R}$ is defined as
\begin{eqnarray}
I& := & \begin{cases}[0,\min_{t}\{x_h(t;x_0)=b\}) & \text{when}~
\{x_h(t;x_0)=b\} \neq \emptyset; \\
[0,\infty) & \text{otherwise}.\end{cases}
\end{eqnarray}
\end{lemma}

\begin{proof}
Under the hypotheses of the lemma we have that $h(x_0) > g(x_0)$,
and it follows that for small enough $t$, $x_h(t)
> x_g(t)$ (omitting the argument $x_0$ for notational
convenience). Note also that $x_h$ and $x_g$ are increasing whenever
$x_h(t;x_0)<b$ and $x_g(t;x_0)<b$ respectively. Now suppose there
exists $t_2>0$ such that $x_g(t_2)
> x_h(t_2)$ and $x_g(t_2)<b$. Then by continuity there exists
$t_1<t_2$ such that $a < x_h(t_1) = x_g(t_1) < b$ and $h(x_h(t_1))
\leq g(x_g(t_1))$. Define $x' := x_h(t_1) = x_g(t_1)$. It follows
that $g(x') \geq h(x')$, which contradicts our starting hypothesis.
We conclude that $x_h(t)\geq x_g(t)$ for all $t\in I$.
\end{proof}

We have the following result:
\begin{corollary}
Let $\theta(\cdot)$ be a solution to the homogeneous
system~(\ref{eq:homsystem}) with initial condition
$\theta(t_0)=\theta_0$. Then
\begin{equation}\label{eq:Rupperbound}
     D(t):=\frac{1}{1-e^{-2k(t-t_0)}\left(\frac{L(\theta_0)-1}{L(\theta_0)}\right)}
\end{equation}
is a dominating function for $L(\theta(t))$ on
$[t_0,\infty)$.
\end{corollary}
\begin{proof}
    By~Lemma~\ref{lem:lowerbound} we have that
    $
        \dot{L}(\theta(t)) \leq 2k
        L(\theta(t))\left(1-L(\theta(t))\right)
    $
    for all~$t$. We claim that, on $[t_0,\infty)$, $L(\theta(t))$ is dominated by the solution $y(t)$ of the ODE
\begin{equation}
    \left\{\begin{array}{ccc}\dot{y}  & =  & 2k y(1-y) \\
    y(t_0) & =  & L(\theta_0)\end{array}\right.
\end{equation}
which is given as
\begin{equation}
y(t) =
\frac{1}{1-e^{-2k(t-t_0)}\left(\frac{L(\theta_0)-1}{L(\theta_0)}\right)},\quad
t\geq t_0.
\end{equation}
To prove this, suppose $\dot{L}(\theta(t)) = 2k
        L(\theta(t))\left(1-L(\theta(t))\right)$ for some $t\geq t_0$ and let $t'$ denote the smallest such $t$ (in case $\dot{L}(\theta(t)) <
2k
        L(\theta(t))\left(1-L(\theta(t))\right)$ for all
        $t\geq t_0$, the result follows immediately from Lemma~\ref{lem:dominator}). Then by
        Lemma~\ref{lem:equalalltimes}, we have that
$
        \dot{L}(\theta(t)) = 2k
        L(\theta(t))\left(1-L(\theta(t))\right)
$
        for all $t\geq t'$, and it follows that
\begin{eqnarray}\label{eq:hello}
\left.L(\theta(t))\right|_{L(\theta(t'))=a}=\frac{1}{1-e^{-2k(t-t')}\left(\frac{a-1}{a}\right)},
\quad t\geq t'.
\end{eqnarray}
One can easily verify that
$\left.L(\theta(t))\right|_{L(\theta(t'))=a}$ is
nondecreasing as a function of $a$ for all $t\geq t'$ (and
$a\in [0,1])$. Now let $l$ be an upper bound for
$L(\theta(t'))$. It follows that
$\frac{1}{1-e^{-2k(t-t')}\left(\frac{l-1}{l}\right)}$ is a
dominating function for $L(\theta(t))$ on the interval
$[t',\infty)$. To compute an upper bound for
$L(\theta(t')$, we proceed as follows. By definition
of~$t'$, we have that $\dot{L}(\theta(t)) < 2k
        L(\theta(t))\left(1-L(\theta(t))\right)$ for all
        $t< t'$. It follows from Lemma~\ref{lem:dominator} that
        \begin{eqnarray}\label{eq:hello2}
L(\theta(t))& \leq &
\frac{1}{1-e^{-2k(t-t_0)}\left(\frac{L(\theta_0)-1}{L(\theta_0)}\right)},
\quad t_0\leq t < t'
\end{eqnarray}
By continuity, we have that \[L(\theta(t'))\leq
\frac{1}{1-e^{-2k(t'-t_0)}\left(\frac{L(\theta_0)-1}{L(\theta_0)}\right)}.\]
Using the upper bound
$\frac{1}{1-e^{-2k(t'-t_0)}\left(\frac{L(\theta_0)-1}{L(\theta_0)}\right)}$
for $L(\theta(t'))$, it follows from~(\ref{eq:hello}) that
\begin{eqnarray}\label{eq:helloagain}
L(\theta(t)) & \leq &
\frac{1}{1-e^{-2k(t-t_0)}\left(\frac{L(\theta_0)-1}{L(\theta_0)}\right)},\quad
t\geq t'
\end{eqnarray}
Combining~(\ref{eq:hello2}) and (\ref{eq:helloagain}), we
arrive at the desired result. This concludes the proof.

\end{proof}

\begin{figure}[thpb]
      \centering
      \psfrag{A}[rc][rc][0.8][0]{$0$}
      \psfrag{B}[rc][rc][0.8][0]{$0.2$}
      \psfrag{C}[rc][rc][0.8][0]{$0.4$}
      \psfrag{D}[rc][rc][0.8][0]{$0.6$}
      \psfrag{E}[rc][rc][0.8][0]{$0.8$}
      \psfrag{F}[rc][rc][0.8][0]{$1$}

      \psfrag{G}[tc][tc][0.8][0]{$0$}
      \psfrag{H}[tc][tc][0.8][0]{$1$}
      \psfrag{I}[tc][tc][0.8][0]{$2$}
      \psfrag{J}[tc][tc][0.8][0]{$3$}
      \psfrag{K}[tc][tc][0.8][0]{$4$}
      \psfrag{L}[tc][tc][0.8][0]{$5$}

      \psfrag{M}[tc][tc][1][0]{$t\rightarrow$}
      \psfrag{N}[lc][lc][1][0]{\hspace{-1.1cm}$\begin{array}{c}\uparrow \\ L(\theta(t))\end{array}$}

      \includegraphics[width=6cm]{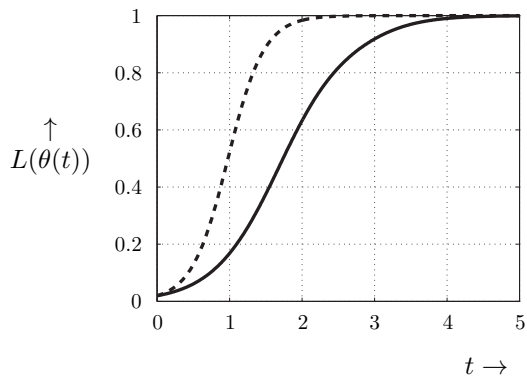}
      \caption{Numerical simulation of the homogeneous system~(\ref{eq:homsystem}) with $N=100$ oscillators and coupling coefficient $k=2$: time evolution of $L(\theta(t)):=R^2(\theta(t))$ (solid line) and the dominating function $D(t)$|Eqn.~(\ref{eq:Rupperbound}) (dashed line).  }
      \label{fig:Rupperbound}
 \end{figure}

Figure~\ref{fig:Rupperbound} shows the graph of
$L(\theta(t))$ and that of the dominating function
$D(t)$|Eqn.~(\ref{eq:Rupperbound}) for a particular
realization of the initial condition $\theta_0$. In this
example, $N=100$ and $k=2$. We observe that, in agreement
with Conjecture~\ref{con:globalstable}, the solution
converges to a globally phase-locked state, that is
$L((\theta(t))\rightarrow 1$. Note that convergence can be
very slow depending on the choice of initial condition.
Indeed, for any $T\in \mathbb{R}$ and any $\epsilon>0$, we
can find $\delta>0$ such that if $L(\theta_0)< \delta$ then
$L(\theta(t))< \epsilon$ for all $t\leq t_0+T$. The upshot
of this is that if the initial condition $\theta_0$ is
selected by drawing from a uniform distribution and the
number of oscillators is large, then $L(\theta_0)$ is
likely to be small, and as a consequence convergence to the
stable equilibrium is likely to be slow. In the limit case
when $N$ tends to infinity, we have that $L(\theta_0)$
tends to zero with probability $1$ and the time required
for $L(\theta(t))$ to exceed some given finite threshold
diverges to infinity.

Let
$\sigma_\omega:=\sqrt{\frac{1}{N}\sum_{j=1}^N\left(\omega_j-\langle
\omega \rangle\right)^2}$ denote the (sample) standard
deviation associated with the vector of intrinsic
frequencies $\omega$. Using Lemma~\ref{lem:lowerbound} we
derive another lower bound on the critical coupling, as
follows:

\begin{corollary}\label{eq:lowerboundl2}
The critical coupling $k_\text{c}$ satisfies
\begin{eqnarray}
    k_\text{c} & \geq & 2\sigma_\omega
\end{eqnarray}
\end{corollary}
\begin{proof}Let $x^*\in V\mathbb{R}^N$ be a
fixed point of the system~(\ref{eq:kuramotomodelgrounded}).
Then by definition $k\|f(x^*)\|_2 = \|V\omega\|_2 =
\sqrt{N}\sigma_\omega$ and by Lemma~\ref{lem:lowerbound} we
have that
\begin{equation}\label{eq:upperboundapplication}
    \|f(x^*)\|_2 \leq
    \sqrt{N}\sqrt{R^2(x^*)\left(1-R^2(x^*)\right)}.
\end{equation}
It is not hard to see that the right hand side of
(\ref{eq:upperboundapplication}) is upper bounded by
$\frac{1}{2}\sqrt{N}$. It follows that
\begin{equation}
    k \geq \frac{\sqrt{N}
    \sigma_\omega}{\sqrt{N}\frac{1}{2}} = 2\sigma_\omega
\end{equation}
This completes the proof.
\end{proof}
Note that Corollary~\ref{eq:lowerboundl2} is in agreement
with the intuition that greater variation in intrinsic
frequencies requires stronger coupling to achieve global
phase-locking.

\begin{figure}[t]
\centering
 \psfrag{A}[ct][ct][0.8][0]{$-0.5$}
\psfrag{B}[ct][ct][0.8][0]{$0.5$}
\psfrag{O}[ct][ct][0.8][0]{$0$}

\psfrag{C}[cr][cr][0.8][0]{$0$}
\psfrag{E}[cr][cr][0.8][0]{$0.5$}
\psfrag{D}[cr][cr][0.8][0]{$1$}

\psfrag{N}[cr][cr][1][0]{$\begin{array}{c}\uparrow
\\R(x^*)\end{array}$}
\psfrag{M}[ct][ct][1][0]{$(\sigma_\omega/{k})\rightarrow$}
\includegraphics[width=6cm]{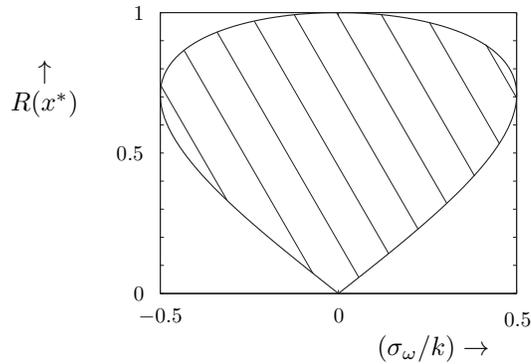}
\caption{Graph associated with inequality~(\ref{eq:Rineq}).
For a given value of the ratio $({\sigma_\omega}/{k})$, the
magnitude of the order parameter $R(\cdot)$, evaluated at a
fixed point $x^*$, must lie within the striped region.}
      \label{fig:tol}
 \end{figure}

Using Lemma~\ref{lem:lowerbound} we can compute bounds on
the value of the order parameter evaluated at the fixed
points of the system, should they exist. Indeed, suppose
$k>k_\text{c}$, then for any fixed point $x^*\in
V\mathbb{R}^N$ we have that
\begin{equation}
    \sqrt{R^2(x^*)\left(1-R^2(x^*)\right)} \geq
    \frac{\sigma_\omega}{k}.
\end{equation}
Solving for $R(x^*)$ gives
\begin{equation}\label{eq:Rineq}
\frac{1}{2}-
\frac{1}{2}\sqrt{1-4\left(\frac{\sigma_\omega}{k}\right)^2}\leq
R^2({x}^*) \leq \frac{1}{2}+
\frac{1}{2}\sqrt{1-4\left(\frac{\sigma_\omega}{k}\right)^2}
\end{equation}
The graph associated with inequality~(\ref{eq:Rineq}) is
shown in Figure~\ref{fig:tol}.

\section{Necessary and Sufficient conditions}
In the last section, we derived lower bounds for the
critical coupling of the system
(\ref{eq:kuramotomodelgrounded}), which provide necessary
conditions for the existence of fixed points. We next
derive conditions that are both necessary and sufficient
for fixed points to exist, and we shall use these results
to describe an algorithm for computing the critical
coupling in Section~6.

Throughout this section, the function $f: \mathbb{R}^N \mapsto
\mathbb{R}R^N$ is given by (\ref{eq:defoffstandard}) and the set
$\mathcal{F}(k,\Omega)$ is defined as
\[
\mathcal{F}(k,\Omega) := \left\{x\in V\mathbb{R}^N: kf(x) =
-\Omega\right\}, \quad k\geq 0,\quad \Omega\in
V\mathbb{R}^N
\]
On $\mathcal{F}(k,\Omega)$ we introduce a notion of
equivalence, as follows.
\begin{defn}[Equivalence on $\mathcal{F}(k,\Omega)$]\label{def:equivalence}Given $\Omega\in V\mathbb{R}^N$ and $k\geq 0$, let $x,x'\in
F(k,\Omega)$.  We say that $x$ and $x'$ are equivalent
($x\simeq x'$) if $R(x) = R(x')$.
\end{defn}
To motivate this definition consider the following fact.
Let $k,\Omega$ be given, and let $s\in \mathbb{Z}^N$ be
such that $\sum_{j=1}^Ns_j = 0$. If $x$ is a fixed point of
the system~(\ref{eq:kuramotomodelgrounded}), then
$x':=x+2s\pi$ is also a fixed point of the
system~(\ref{eq:kuramotomodelgrounded}) and, in addition,
$R(x') = R(x)$.

The following theorem provides a necessary and sufficient
condition for the system~(\ref{eq:kuramotomodelgrounded})
to have a fixed point, given a particular coupling strength
$k$, and a particular realization of intrinsic frequencies,
$\Omega$.

\begin{theorem}\label{thm:necsuf}
Let $k>0$ and $\Omega\in V\mathbb{R}^N$. Then
$\mathcal{F}(k,\Omega)\neq \emptyset$ iff there
    exists $\beta \in [\frac{1}{k}\|\Omega\|_\infty,1]\subset\mathbb{R}$ and $a\in \{-1,1\}^N$ such that
    \begin{equation}\label{eq:betacondition}
     \beta=\frac{1}{N}\sum_{j=1}^Na_j\sqrt{1-\left(\frac{\Omega_j}{k \beta}\right)^2}.
     \end{equation}
     Moreover, suppose $(a^1,\beta^1)$ and $(a^2,\beta^2)$ both satisfy~(\ref{eq:betacondition}) and let $x^1, x^2\in
     F(k,\Omega)$ be such that
     \begin{equation}\label{eq:geenzin}\left\{\begin{array}{rcc}k \beta^i\sin(\psi(x^i) - x^i_j) & = &
     -\Omega_j\\
     a_j \cos(\psi(x^i) - x^i_j) & \geq &
     0\end{array}\right.,\quad i\in \{1,2\},\quad j=1,2,\ldots N
     \end{equation}
     Then $x^1\simeq x^2$ iff $\beta^1 = \beta^2$ and
     $\sum_{j=1}^N\left(a_i-a_j\right)\sqrt{1-\left(\frac{\Omega_j}{k
     \beta^1}\right)^2} = 0$.
\end{theorem}

\begin{proof}
Suppose $\Omega\neq 0$ (the case $\Omega=0$ is easy). Let
${x}^*\in V\mathbb{R}^N$ be a fixed point of
(\ref{eq:kuramotomodelgrounded}). By definition,
$kf({x}^*)=-\Omega$, and since $\Omega\neq 0$, we have that
$f({x}^*)\neq 0$, and consequently $R({x}^*)\neq 0$. It follows
that
     \begin{equation}\label{eq:fixedpointidentity}
        \sin(\psi({x}^*)-{x}^*_i) = -\frac{\Omega_i}{k
        R({x}^*)}, \quad\quad i=1,2,\ldots,N.
     \end{equation}
     Let $\beta:=R({x}^*)$. By~(\ref{eq:fixedpointidentity}) we have that $\beta\geq \frac{1}{k}\|\Omega\|_\infty$.
     Recall that for all $x\in \mathbb{R}^N\backslash
     \mathcal{R}_0$, $R(x)$ can be written as
     \begin{equation}\label{eq:Ridentity}
     R(x)=\frac{1}{N}\sum_{j=1}^N\cos(\psi(x)-x_j),
     \end{equation} and let $a_i$ be given as
     \begin{equation}\label{eq:sign}
     a_i:=\begin{cases} -1 & \text{if $\cos(\psi({x}^*)-{x}^*_i)
     \leq
     0$;} \\
     +1 & \text{otherwise.}
     \end{cases}
     \end{equation}
     Combining~(\ref{eq:fixedpointidentity}), (\ref{eq:Ridentity}), and (\ref{eq:sign}), we
     arrive at
     \begin{equation}\label{eq:condition}
     \beta=\frac{1}{N}\sum_{j=1}^Na_j\sqrt{1-\left(\frac{\Omega_j}{k \beta}\right)^2}.
     \end{equation}
     This proves necessity. To prove sufficiency, let $a\in \{-1,1\}^N$ be given, and suppose $\beta\geq
     \frac{1}{k}\|\Omega\|_\infty>0$ (again, the case $\Omega=0$ is easy). Then for every $c\in \mathbb{R}$, the system
\begin{equation}\label{eq:auxsystem}
    \left\{\begin{array}{ccc}k\beta\sin(-y_i-c) & = &
    -\Omega_i\\ a_i\cos(-y_i-c) & \geq  & 0,
    \end{array}\right.\quad\quad i=1,2,\ldots,N
\end{equation}
has a unique solution ${y}^*\in [-\pi,\pi)^N$. We pick $c$
such that $\sum_{j=1}^N y^*_j = 0$. Since
$\sum_{j=1}^N\sin(y^*_j+c)=0$, it follows that
\begin{eqnarray}\label{eq:interme}
R(y^*) = R(y^*+c{\bf 1}) =
\left|\sum_{j=1}^N\cos(y^*_j+c)\right|
\end{eqnarray}
From~(\ref{eq:auxsystem}), we have that
\begin{eqnarray}\label{eq:interme2}
\cos(y^*_i+c) = a_i\sqrt{1 -
\left(\frac{\Omega_i}{k\beta}\right)^2}\quad i=1,\ldots,N.
\end{eqnarray}
Combining~(\ref{eq:interme}) and (\ref{eq:interme2}), we
arrive at
\begin{equation}
    R({y}^*) =
    \left|\frac{1}{N}\sum_{j=1}^Na_j\sqrt{1-\left(\frac{\Omega_j}{k\beta}\right)^2}\right|
\end{equation}
The second part of the theorem follows easily after noting
that if $x^i$ satisfies~(\ref{eq:geenzin}) then $R(x^i) =
\beta^i$, $i=1,2$.
\end{proof}

Theorem~\ref{thm:necsuf} gives us a necessary and sufficient
condition for the equation $kf(x)=-\Omega$ to have at least one
solution for a given value of $k$. It is not clear, however, that
there exists a $k$ for which this condition is satisfied. The
following Corollary provides an easy sufficient condition.
\begin{corollary}\label{cor:sufficient}
    Let $k> 0$ and  $\Omega \in V\mathbb{R}^N$. Suppose
\begin{equation}\label{eq:conditionsuf}
    \frac{1}{k}\|\Omega\|_\infty \leq
    \frac{1}{N}\sum_{j=1}^Na_j\sqrt{1-\left(\frac{\Omega_j}{\|\Omega\|_\infty}\right)^2}.
\end{equation}
for some $a\in \{-1,1\}^N$. Then $\mathcal{F}(k,\Omega)\neq
\emptyset$.

\end{corollary}
\begin{proof}
Suppose $\Omega\neq 0$ (again, the case $\Omega=0$ is
easy). Let $a\in \{-1,1\}^N$. Define
$m:[\frac{1}{k}\|\Omega\|_\infty,1]\mapsto \mathbb{R}$,
$m(\beta):=\beta$ and
$n:[\frac{1}{k}\|\Omega\|_\infty,1]\times \{-1,1\}^N
\mapsto \mathbb{R}$,
\begin{equation}
    n(\beta,a):=\frac{1}{N}\sum_{j=1}^Na_j\sqrt{1-\left(\frac{\Omega_j}{k
\beta}\right)^2}.
\end{equation}
Since $\Omega\neq 0$ we have that $m(1)>n(1,a)$. Now
suppose condition~(\ref{eq:conditionsuf}) is satisfied.
Then we have that $m(\frac{1}{k}\|\Omega\|_\infty)\leq
n(\frac{1}{k}\|\Omega\|_\infty,a)$, and by the Intermediate
Value Theorem there must exist $\beta^*\in
[\frac{1}{k}\|\Omega\|_\infty,1]$ such that $m(\beta^*) =
n(\beta^*,a)$. It follows from Theorem~\ref{thm:necsuf}
that the system~(\ref{eq:kuramotomodelgrounded}) has a
fixed point.
\end{proof}

\begin{corollary}\label{cor:numberoffixedpoints}
Let $\Omega\in V\mathbb{R}^N$. Then: $(i)$ the critical
coupling $k_\text{c}$ is finite; $(ii)$ for large enough
coupling, the system~(\ref{eq:kuramotomodelgrounded}) has
at least $2^{N-1}$ fixed points.
\end{corollary}
\begin{proof}
Note that the right hand side of (\ref{eq:conditionsuf})
does not depend on $k$. Hence, it follows that, provided
\begin{equation}\label{eq:positivenegative}
\frac{1}{N}\sum_{j=1}^Na_j\sqrt{1-\left(\frac{\Omega_j}{\|\Omega\|_\infty}\right)^2}
> 0,
\end{equation}Condition~(\ref{eq:conditionsuf}) is always satisfied
for large enough $k$. Furthermore, it follows easily that
if (\ref{eq:positivenegative}) is not satisfied for given
$a$, then it is satisfied for $a':=-a$. This implies ($i$)
that the critical coupling $k_\text{c}$ is always finite,
and ($ii$) that the set $A^{+}:=\{a\in
\{-1,1\}^N:\text{Eqn.~(\ref{eq:positivenegative}) is
satisfied}\}$ contains precisely $2^{N-1}$ elements
(counting multiplicity), each of which defines a unique (up
to equivalence in the sense of
Definition~\ref{def:equivalence}) fixed point. This
concludes the proof.\end{proof}

\begin{corollary}\label{cor:synchronized}
Let $k>0$ and $\Omega \in V\mathbb{R}^N$. Then
$\mathcal{F}(k,\Omega)\neq \emptyset$ if and only if there
exist $\beta \in [\frac{1}{k}\|\Omega\|_\infty,1]$ such
that
\[
\beta=\frac{1}{N}\sum_{j=1}^N\sqrt{1-\left(\frac{\Omega_j}{k
\beta}\right)^2}.\]
\end{corollary}
\begin{proof}
The proof of Corollary~\ref{cor:sufficient} suggests that
if the fixed point equation~(\ref{eq:betacondition}) does
not have a solution, then necessarily
\[\beta > \frac{1}{N}\sum_{j=1}^Na_j\sqrt{1-\left(\frac{\Omega_j}{k
\beta}\right)^2}\] for all $\beta \in
[\frac{1}{k}\|\Omega\|_\infty,1]$ and for all $a\in
\{-1,1\}^N$. Since we have that
\[\frac{1}{N}\sum_{j=1}^N\sqrt{1-\left(\frac{\Omega_j}{k
\beta}\right)^2}\geq
\frac{1}{N}\sum_{j=1}^Na_j\sqrt{1-\left(\frac{\Omega_j}{k
\beta}\right)^2}\] for all $a\in \{-1,1\}^N$, it follows
that the given condition is necessary and sufficient for
the system~(\ref{eq:kuramotomodelgrounded}) to have at
least one fixed point. This concludes the proof.
\end{proof}

The next and final corollary gives us an upper bound on the
critical coupling.

\begin{corollary}\label{cor:upperbound}
    The critical coupling, $k_\text{c}$, satisfies:

    \begin{equation}
        k_\text{c} \leq
        \frac{\|\Omega\|_\infty}{\frac{1}{N}\sum_{j=1}^N\sqrt{1-\left(\frac{\Omega_j}{\|\Omega\|_\infty}\right)^2}}.
    \end{equation}
\end{corollary}
\begin{proof}
Follows directly from Corollary~\ref{cor:sufficient}.
\end{proof}

\section{An algorithm for computing $k_\text{c}$}
In this section we present a bisection algorithm that will
allow us to numerically evaluate the critical coupling with
arbitrary precision. Throughout, we shall assume that
$\Omega\neq 0$. Define
$\mathcal{I}:=(\|\Omega\|_{\infty},\infty)$, and let
$p_i:\mathcal{I}\mapsto (0,1]$, and $P:\mathcal{I}\mapsto
(0,1]$ be given as
\begin{equation}\label{eq:CapitalP}
p_i(u):=\sqrt{1-\left(\frac{\Omega_i}{u}\right)^2}\quad;\quad
P(u):=\frac{1}{N}\sum_{j=1}^N p_j(u)\end{equation}Also,
define $h(u;k):\mathcal{I}\times \mathbb{R}_{+}\mapsto
\mathbb{R}_{+}$,
\begin{equation}\label{eq:straightline}
h(u;k):=\frac{1}{k}u.
\end{equation}
From Corollary~\ref{cor:synchronized} it follows that the
critical coupling is the smallest $k$ for which the
equation $P(u) = h(u;k)$ has at least one solution on
$\mathcal{I}$. We have the following result.

\begin{theorem}\label{thm:characterizesolution}
For all $\Omega\in V\mathbb{R}^N$, $\Omega\neq 0$, the
equation
\begin{equation}\label{eq:keyequation}
2\frac{1}{N}\sum_{j=1}^N\sqrt{1-\left(\frac{\Omega_j}{u}\right)^2}
=
\frac{1}{N}\sum_{j=1}^N\frac{1}{\sqrt{1-\left(\frac{\Omega_j}{u}\right)^2}}.
\end{equation}
has a unique solution $u^*\in \mathcal{I}$, and we have
that
\begin{equation}\label{eq:defkc}
k_\text{c}=\frac{u^*}{\frac{1}{N}\sum_{j=1}^N\sqrt{1-\left(\frac{\Omega_j}{u^*}\right)^2}.
}
\end{equation}
\end{theorem}
\begin{proof}
Observe that, by strict concavity of $P$ and linearity of
$h(\:\cdot\:;k)$, the equation $P(u) = h(u;k)$ can have at
most two solutions on $\mathcal{I}$ for any $k>0$. We shall
now show that, when $k=k_\text{c}$, it can have no more
than one solution. Since, by definition of critical
coupling, $P(u) = h(u;k_\text{c})$ must have \emph{at
least} one solution, we shall conclude that it has
precisely one solution. Let $k=k_\text{c}$ and suppose
there exist $u^1, u^2\in \mathcal{I}$, $u^1\neq u^2$, such
that $P(u^1) = h(u^1;k_\text{c})$ and
$P(u^2)=h(u^2;k_\text{c})$. By strict concavity of $P$ we
have that $P(\frac12(u^1+u^2))>
\frac12\left(P(u^1)+P(u^2)\right)$. Define $u' :=
\frac12(u^1+u^2)$ and note that $u'\in \mathcal{I}$. We
have that $P(u')> h(u';k_\text{c})$. This implies that
there exists $k'< k_\text{c}$ such that $P(u') = h(u';k')$.
But by definition $k_\text{c}$ is the smallest $k$ for
which $P(u) = h(u;k)$ has a solution. We arrive at a
contradiction and conclude that $u^1=u^2$; or in other
words, that the equation $P(u) = h(u;k_\text{c})$ has
exactly one solution on $\mathcal{I}$. Denoting this
solution by $u^*$, it is not hard to see that, at $u=u^*$,
the derivative of $P$ with respect to $u$ and the
derivative of $h$ with respect to $u$ (both of which are
defined on the entire interval $\mathcal{I}$) must
coincide. For suppose $\frac{\partial h}{\partial u}(u^*) <
\frac{\partial P}{\partial u}(u^*)$, then by continuity
there exists $\delta>0$ such that $h(u;k_\text{c}) < P(u)$
for all $u$ such that $u-u^*<\delta$. Let $u'$ be one such
$u$. It follows that there exists $k'<k_\text{c}$ such that
$P(u') = h(u';k')$. This leads to a contradiction and we
conclude that $h(u;k_\text{c}) \geq P(u)$. By analogy we
have that $h(u;k_\text{c}) \leq P(u)$. We conclude that
$\frac{\partial h}{\partial u}(u^*) = \frac{\partial
P}{\partial u}(u^*)$. That is,
\begin{eqnarray}
    \frac{1}{{k}_\text{c}} & = &
    \frac{1}{u^*}\frac{1}{N}\sum_{j=1}^N\frac{\left(\frac{\Omega_j}{u^*}\right)^2}{\sqrt{1-\left(\frac{\Omega_j}{u^*}\right)^2}}.
\end{eqnarray}
Or equivalently,
\begin{eqnarray}\label{eq:eq1}
    \frac{u^*}{{k}_\text{c}} & = &
    \frac{1}{N}\sum_{j=1}^N\frac{\left(\frac{\Omega_j}{u^*}\right)^2}{\sqrt{1-\left(\frac{\Omega_j}{u^*}\right)^2}}\nonumber\\
    & = &
    -
    \frac{1}{N}\sum_{j=1}^N\sqrt{1-\left(\frac{\Omega_j}{u^*}\right)^2}+\frac{1}{N}\sum_{j=1}^N\frac{1}{\sqrt{1-\left(\frac{\Omega_j}{u^*}\right)^2}}
\end{eqnarray}
Now recall that by definition of $u^*$, we have that
\begin{eqnarray}\label{eq:eq2}
\frac{u^*}{{k}_\text{c}} & = &
\frac{1}{N}\sum_{j=1}^N\sqrt{1-\left(\frac{\Omega_j}{u^*}\right)^2}.
\end{eqnarray}
Equating the right hand side of Eqn.~(\ref{eq:eq1}) with
the right hand side of Eqn.~(\ref{eq:eq2}) gives
\begin{equation}
2\frac{1}{N}\sum_{j=1}^N\sqrt{1-\left(\frac{\Omega_j}{u^*}\right)^2}
=
\frac{1}{N}\sum_{j=1}^N\frac{1}{\sqrt{1-\left(\frac{\Omega_j}{u^*}\right)^2}}.
\end{equation}
This shows that $u^*$ is a solution
to~(\ref{eq:keyequation}). What remains to be shown is that
$u^*$ is the only solution. Define $v,w: \mathcal{I}\mapsto
\mathbb{R}$,
\[
v(u):=2\frac{1}{N}\sum_{j=1}^N\sqrt{1-\left(\frac{\Omega_j}{u}\right)^2},\quad
w(u):=\frac{1}{N}\sum_{j=1}^N\frac{1}{\sqrt{1-\left(\frac{\Omega_j}{u}\right)^2}},\]
and note that, on their respective domains, $v$ is strictly
monotonically decreasing while $w$ is strictly
monotonically increasing. In addition, note that there
exist $a,b\in \mathcal{I}$ such that $v(a)>w(a)$ and
$v(b)<w(b)$. Hence, by continuity, there must exist a point
$u'\in (a,b)\subset \mathcal{I}$ such that $v(u')=w(u')$.
Monotonicity of $v$ and $w$ implies that this point is
unique. It follows that $u^*$ is the unique solution
of~(\ref{eq:keyequation}) on $\mathcal{I}$. And
by~(\ref{eq:eq2}) we have that
\begin{equation}
k_\text{c}=\frac{u^*}{\frac{1}{N}\sum_{j=1}^N\sqrt{1-\left(\frac{\Omega_j}{u^*}\right)^2}
}
\end{equation}
This concludes the proof.
\end{proof}

Based on the result of
Theorem~\ref{thm:characterizesolution}, we define the map
$K:V\mathbb{R}^N\backslash\{0\} \mapsto \mathbb{R}_{+}$,
\begin{equation}
K(\Omega) =
\frac{u^*}{\frac{1}{N}\sum_{j=1}^N\sqrt{1-\left(\frac{\Omega_j}{u^*}\right)^2}
},
\end{equation}
where, as before, $u^*$ denotes the unique solution
of~(\ref{eq:keyequation}) on $\mathcal{I}$, given $\Omega$.
Note that, given any realization of $\omega$ such that
$V\omega\neq 0$, we have that $k_\text{c}=K(V\omega)$. We
have the following corollary.

\begin{corollary}{~}\label{cor:tightbounds}

\begin{enumerate}
\item For all $\Omega\in V\mathbb{R}^N$, $\Omega\neq 0$, we
have that $\|\Omega\|_\infty \leq K(\Omega) \leq
2\|\Omega\|_\infty$;\item there exists $\Omega\in
V\mathbb{R}^N$ such that $K(\Omega) = 2\|\Omega\|_\infty$
if and only if $N$ is even; \item for every $\epsilon>0$
there exist an positive integer $N$ and $\Omega\in
V\mathbb{R}^N$ such that
$|K(\Omega)-\|\Omega\|_\infty|<\epsilon$.
\end{enumerate}
\end{corollary}

\begin{proof}(Part 1). We show that for all $\Omega \neq 0$, the solution
$u^*$ of equation~(\ref{eq:keyequation}) satisfies $u^*
\leq \sqrt{2}\|\Omega\|_\infty$. The result then follows
easily. Let $u':= \sqrt{2}\|\Omega\|_\infty$. Then we have
that
\[2\frac{1}{N}\sum_{j=1}^N\sqrt{1-\left(\frac{\Omega_j}{u}\right)^2}
 > \frac{1}{N}\sum_{j=1}^N\frac{1}{\sqrt{1-\left(\frac{\Omega_j}{u}\right)^2}} \quad \text{for all~}
 u>
 u'\]
It follows that
\[K(\Omega)\leq
\frac{u'}{\frac{1}{N}\sum_{j=1}^N\sqrt{1-\left(\frac{\Omega_j}{u'}\right)^2}}\leq
\frac{\sqrt{2}\|\Omega\|_\infty}{\frac12\sqrt{2}} =
{2}\|\Omega\|_\infty
\]
for all $\Omega \in V\mathbb{R}^N$. The lower bound
$K(\Omega)\geq \|\Omega\|_\infty$ was obtained earlier in
Section~\ref{sec:lowerbounds}. \\

\noindent (Part 2.) From the above it follows that
$K(\Omega)=2\|\Omega\|_\infty$ if and only if
\[2\frac{1}{N}\sum_{j=1}^N\sqrt{1-\left(\frac{\Omega_j}{u'}\right)^2}
 =
 \frac{1}{N}\sum_{j=1}^N\frac{1}{\sqrt{1-\left(\frac{\Omega_j}{u'}\right)^2}}
\]
or equivalently, $\Omega^2_i=\Omega^2_j$ for all $(i,j)$.
It is easy to see that this latter condition is never
satisfied when $N$ is odd (keeping in mind that
$\sum_j\Omega_j=0$). Now suppose $N$ is even and pick any
$c\neq 0$. Define
\[ \Omega_i :=
\begin{cases}
c& i=1,2,\ldots \frac{N}{2};\\
-c & i=\frac{N}{2}+1,\ldots, N
\end{cases}.
\] Then we have that $\Omega\in V\mathbb{R}^N$. Moreover,
$\Omega^2_i=\Omega^2_j=c^2$ for all $(i,j)$. It follows
that $K(\Omega)=2\|\Omega\|_\infty$.

\noindent(Part 3.) Let $\epsilon > 0$ be given and suppose
$N$ is odd. Pick $c\neq 0$ and define
\[ \Omega_i :=
\begin{cases}
0& i=1,2,\ldots,N-1;\\
c & i = N.
\end{cases}
\]
Then~(\ref{eq:keyequation}) evaluates to
\begin{equation}\label{eq:keystep}
2(N-1)+ 2\sqrt{1-\left(\frac{c}{u}\right)^2} =
(N-1)+\frac{1}{\sqrt{1-\left(\frac{c}{u}\right)^2}}.
\end{equation}
and it is not hard to see that as $N$ tends to infinity,
the solution $u^*$ of~(\ref{eq:keystep}) tends to $c$.
Indeed, for $N\geq 2$ we have
\begin{eqnarray}
\left(\frac{c}{u^*}\right)^2 & = & \frac12
(N-1)\left(-\frac14 (N-1)+ \frac14\sqrt{(N-1)^2+8}\right).
\end{eqnarray}
Let $\epsilon_1>0$, and pick $N$ such that
$\frac{1}{N}<\epsilon_1$ and $u^*< (1+\epsilon_1)c$. It
follows that
\begin{eqnarray}
\frac{1}{N}\sum_{j=1}^N\sqrt{1-\left(\frac{\Omega_j}{u^*}\right)^2}
& < & 1-\epsilon_1,
\end{eqnarray}
and hence
\begin{equation}\label{eq:smallerthan}K(\Omega)~:=~\frac{u^*}{\frac{1}{N}\sum_{j=1}^N\sqrt{1-\left(\frac{\Omega_j}{u^*}\right)^2}}~<~c\left(\frac{1+\epsilon_1}{1 -
\epsilon_1}\right)~=~c +
2\left(\frac{\epsilon_1}{1-\epsilon_1}\right)c.
\end{equation}
Now let $\epsilon_1$ be given as
$\epsilon_1:=\frac{\epsilon}{2c+\epsilon}$ and choose $N$
accordingly. It follows that $K(\Omega)< \|\Omega\|_\infty
+ \epsilon$. This concludes the proof.

We are now ready to present our algorithm, which, given
$\Omega$, will compute $u^*$ with user-defined precision
$\epsilon>0$ in a finite number of iterations,
$n=\lceil\log_2(\frac{\|\Omega\|_\infty}{\epsilon})\rceil
+1$.
\begin{algorithm}{~}\label{alg:algorithm}

\begin{enumerate}
\item \hspace{0.1cm} $a:=\|\Omega\|_\infty$,
$b:=\sqrt{2}\|\Omega\|_\infty$ \item \hspace{0.1cm} While
$(b-a) > \epsilon$, \item \hspace{0.6cm} $u:=\frac12(b-a)$.
\item \hspace{0.6cm} If $\left[\sum_{j}
\sqrt{\left(1-\frac{\Omega_j}{u}\right)^2}~
>~\frac12 \sum_{j}\frac{1}{
\sqrt{\left(1-\frac{\Omega_j}{u}\right)^2}}\right]$~then
$a:=u$, else $b:=u$.\item \hspace{0.1cm} End.
\end{enumerate}
\end{algorithm}

\end{proof}

Once we have an estimate $\hat{u}$ of $u^*$, we can
use~(\ref{eq:defkc}), replacing $u^*$ with $\hat{u}$, to
estimate $k_\text{c}$.

%
%
%

\section{Numerical Example}
\begin{figure}[thpb]
\centering \psfrag{C}[rc][rc][1][0]{$-2$}
      \psfrag{D}[rc][rc][1][0]{}
      \psfrag{E}[rc][rc][1][0]{}
      \psfrag{F}[rc][rc][1][0]{}
      \psfrag{G}[rc][rc][1][0]{$0$}
      \psfrag{H}[rc][rc][1][0]{}
      \psfrag{I}[rc][rc][1][0]{}
      \psfrag{J}[rc][rc][1][0]{}
      \psfrag{K}[rc][rc][1][0]{}
      \psfrag{L}[rc][rc][1][0]{$1.5$}

      \psfrag{A}[tc][tc][1][0]{$1$}
      \psfrag{B}[tc][tc][1][0]{$20$}

      \psfrag{M}[tc][tc][1][0]{$i\rightarrow$}
      \psfrag{N}[lc][lc][1][0]{$\begin{array}{c}\uparrow\\\Omega^{20}_i\end{array}$}
\begin{tabular}{cc}
      \includegraphics[width=5.5cm]{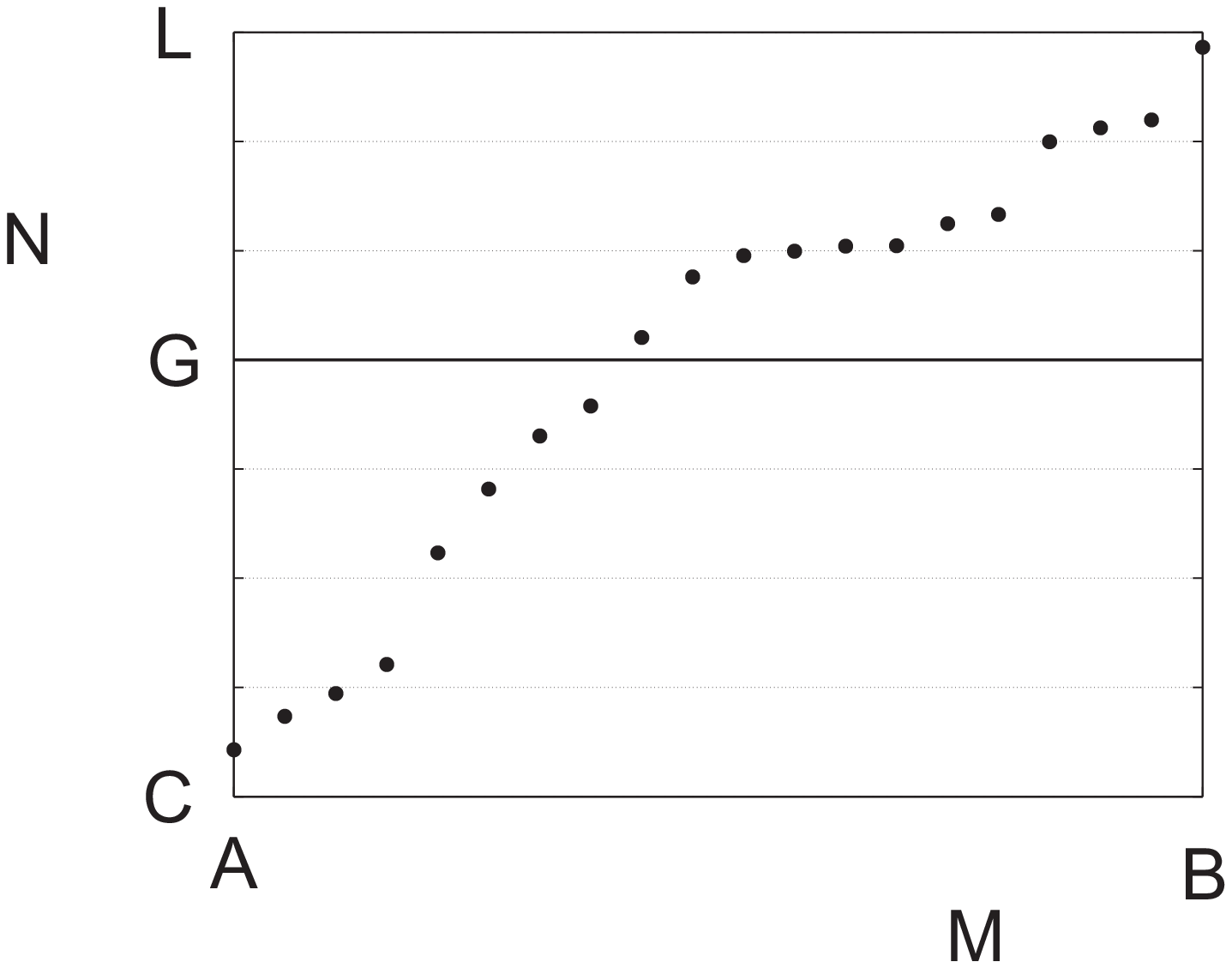}
      &
      \psfrag{B}[tc][tc][1][0]{$200$}
      \psfrag{C}[rc][rc][1][0]{$-2.5$}
       \psfrag{L}[rc][rc][1][0]{$2.5$}
      \psfrag{N}[lc][lc][1][0]{$\begin{array}{c}\uparrow\\\Omega^{200}_i\end{array}$}
      \includegraphics[width=5.5cm]{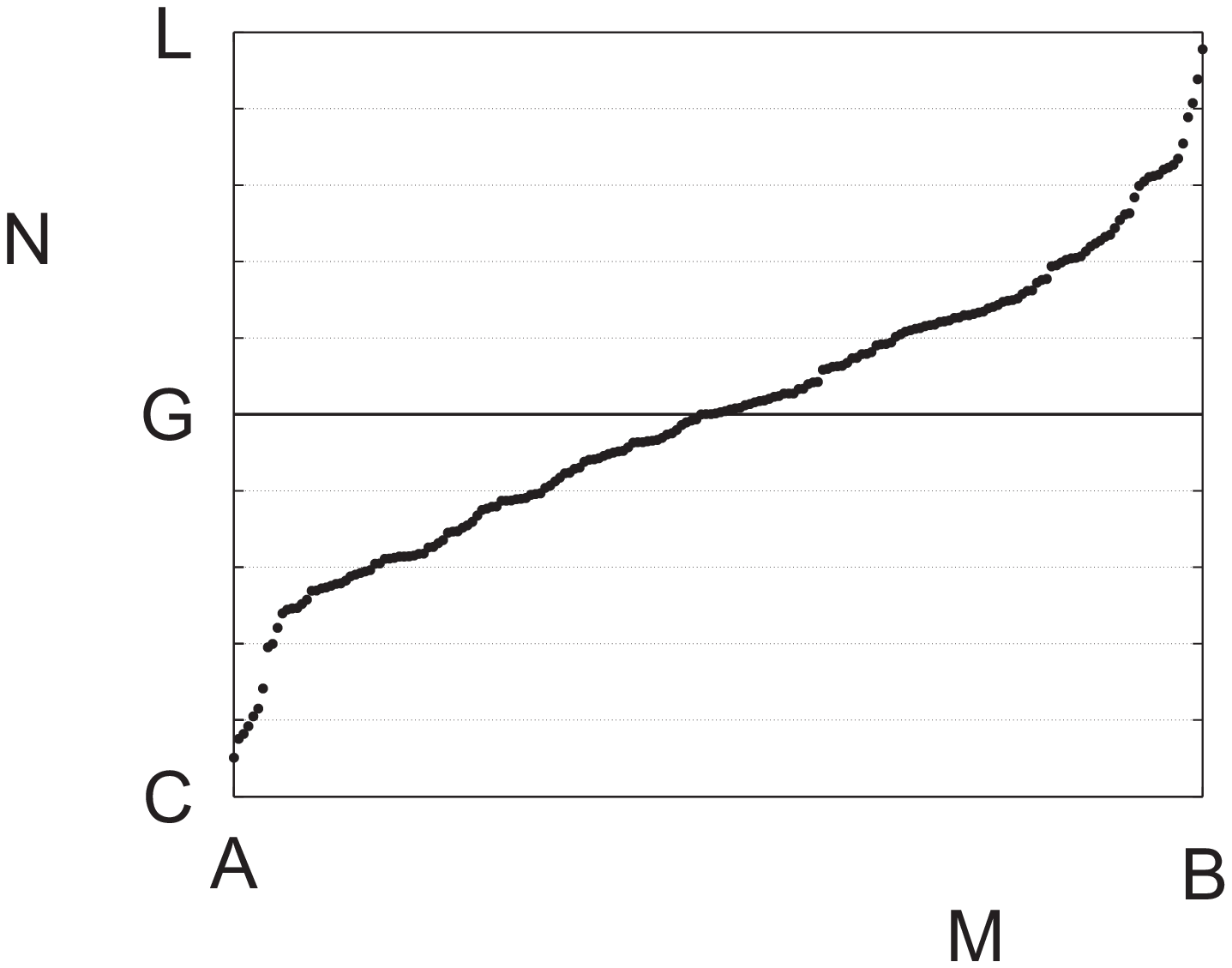}\\
      (a) $N=20$ & (b) $N=200$\end{tabular}
      \caption{The vector of frequencies $\Omega^{N}_i:=\omega^{N}_i-\langle
      \omega^{N}\rangle$, $N\in \{20,200\}$, used in this example. The natural frequencies $\omega^{N}_i$ were sampled from a normal distribution with zero mean and unit variance and relabelled in such a way that $\omega^{N}_1\leq \omega^{N}_2\leq \ldots \leq \omega^{N}_N$.}
      \label{fig:naturalfrequencies}
 \end{figure}

\begin{figure}[thpb]
\centering

\psfrag{M}[ct][ct][0.8][0]{$\beta\rightarrow$}
\psfrag{G}[cr][cr][0.8][0]{$h(k\beta;k)$}
\psfrag{H}[tr][tr][0.8][0]{$P^{20}(k\beta)$}

    \begin{tabular}{cc}

      {\psfrag{A}[ct][ct][0.8][0]{$\frac{1}{k}\|\Omega^{20}\|_\infty$}
\psfrag{B}[ct][ct][0.8][0]{$1$}

\psfrag{C}[cr][cr][0.8][0]{$0.88$}
\psfrag{D}[cr][cr][0.8][0]{$1$}\includegraphics[width=6cm]{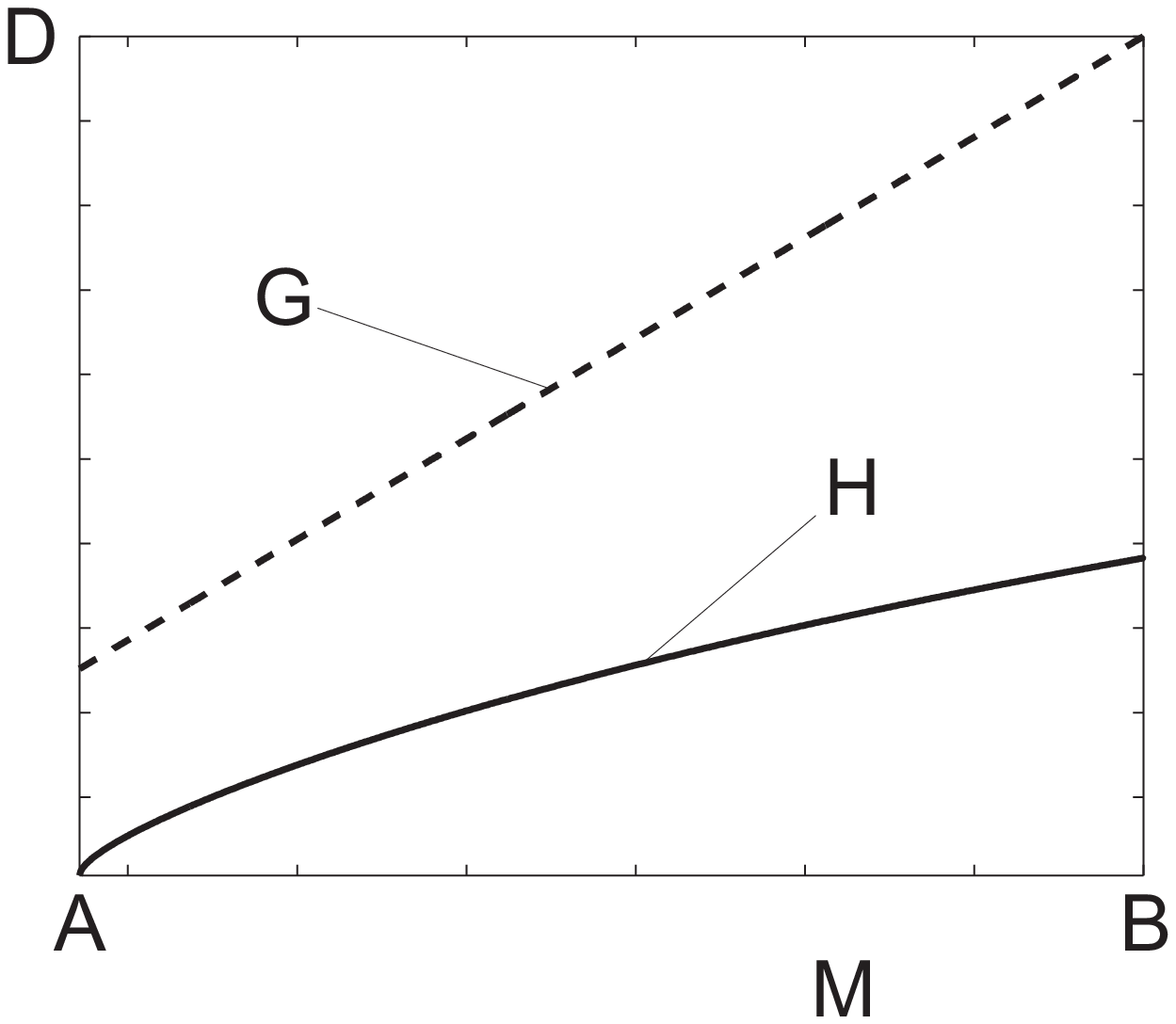}}
      &
      {\psfrag{A}[ct][ct][0.8][0]{$\frac{1}{k}\|\Omega^{20}\|_\infty$}
\psfrag{B}[ct][ct][0.8][0]{$1$}

\psfrag{C}[cr][cr][0.8][0]{$0.88$}
\psfrag{D}[cr][cr][0.8][0]{$1$}\includegraphics[width=6cm]{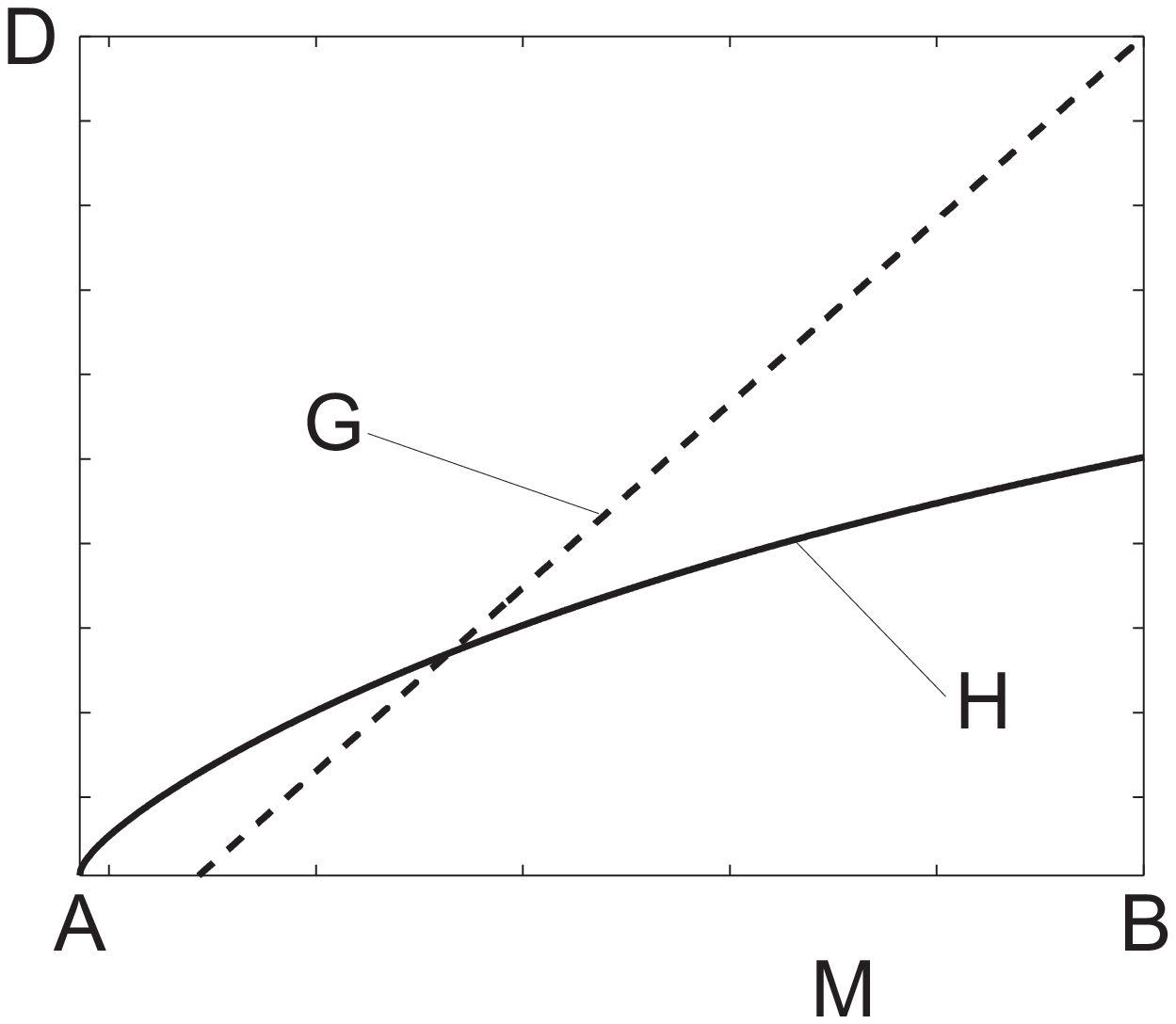}}\\
      (a) $k=2.1$ & (b) $k=2.3$

    \end{tabular}

      \caption{Case $N=20$:~the graph of $P^{20}(k\beta)$ (Eqn.~(\ref{eq:CapitalP})) vs.~$\beta$ for $k=2.1$, $2.3$ and $\beta \in [\frac{1}{k}\|\Omega^{20}\|_\infty,1]$. The dashed line is the graph of $h(k\beta;k)=\beta$ (Eqn.~(\ref{eq:straightline})). An intersection corresponds to a solution of the fixed point equation $P(k\beta)=h(k\beta;k)$, and thus, by Theorem~\ref{thm:necsuf}, to a fixed point of the system~(\ref{eq:kuramotomodelgrounded}).}
      \label{fig:fixedpoint} \end{figure}

We illustrate the results presented in this paper by means
of a numerical example. We consider two systems with $N=20$
and $N=200$ oscillators respectively, with frequencies
$\{\Omega_i^N\}$ as depicted in
Figure~\ref{fig:naturalfrequencies}.
The frequencies in this example were sampled from a normal
distribution with zero mean and unit variance and
relabelled such that $\omega^N_1\leq \omega^N_2 \ldots \leq
\omega^N_N$ (note that this can be done without loss of
generality). For this particular realization of
$\omega^{20}$ ($\omega^{200}$), we have that
$\|\Omega^{20}\|_\infty = 1.7858$ ($\|\Omega^{200}\|_\infty
= 2.3893)$ and
\[\frac{1}{N}\sum_{j=1}^N\sqrt{1-\left(\frac{\Omega_j}{\|\Omega\|_\infty}\right)^2}
= 0.8015~(0.9139).
\]
It follows from Corollary~\ref{cor:upperbound} that
$k_\text{c}\leq 2.2281$ ($2.6145$) and
by~(\ref{eq:lowerboundsimple}), we have that
$k_\text{c}\geq \|\Omega\|_{\infty}=1.7858$ ($2.3893$).
Figure~\ref{fig:Example} shows the time evolution of the
magnitude squared of the order parameter, $R^2(t)$
(previously denoted as $L(t)$), for two different initial
conditions and two values of the coupling coefficient,
$k=2.3$ and $k=2.65$ ($k=2.1$ and $k=2.3$). We observe that
when $k$ is slightly greater than the known lower bound on
$k_\text{c}$, the value of $R^2(t)$ converges to a constant
and inspection shows that the solution $x(t)$ of the
system~(\ref{eq:kuramotomodelgrounded}) tends to a fixed
point. On the other hand, when the coupling coefficient is
slightly below the known upper bound on the critical
coupling, the trajectories $x(t)$ appear not to converge.
Note that in this case we do not know whether the
system~(\ref{eq:kuramotomodelgrounded}) has a fixed point
or not, as the condition stated in
Corollary~\ref{cor:upperbound} is only sufficient while at
the same time the respective coupling strengths exceed
their known lower bounds ($1.7858$ and $2.3893$
respectively). To gain more insight into this situation let
us consider the case $N=20$ in some more detail. We fix the
coupling coefficient at $k=2.1$, and numerically evaluate
the function $P^{20}(k,\cdot)$,
\begin{equation}\label{eq:gbeta}
P^{20}(k\beta)=\frac{1}{20}\sum_{j=1}^{20}\sqrt{1-\left(\frac{\Omega^{20}_j}{k
\beta}\right)^2},
\end{equation}
for several values of $\beta$ in the interval
$[\frac{1}{k}\|\Omega^{20}\|_\infty,1]$. We repeat the same
computation for $k=2.3$. The result is shown in
Figure~\ref{fig:fixedpoint}. We observe that the equation
$P^{20}(k\beta) = \beta$ does not have a solution on the
interval $[\frac{1}{k}\|\Omega^{20}\|_\infty,1]$ when
$k=2.1$, but does have a solution when $k=2.3$.

We use Algorithm~\ref{alg:algorithm} to compute the `exact'
value of the critical coupling to the fifth significant
digit. We find that $k_\text{c}=2.2198$ for the case $N=20$
and $k_\text{c}=2.6144$ for the case $N=200$. Note that in
both cases, but particularly the latter, the upper bounds
($2.2281$ and $2.6145$ respectively) provide good estimates
of the true values of the critical coupling.

\begin{figure}[t]
\centering \psfrag{A}[ct][ct][0.8][0]{$0$}
\psfrag{B}[ct][ct][0.8][0]{$5$}
\psfrag{C}[ct][ct][0.8][0]{$10$}
\psfrag{D}[ct][ct][0.8][0]{$15$}
\psfrag{E}[ct][ct][0.8][0]{$20$}

\psfrag{F}[cr][cr][0.8][0]{$0$}
\psfrag{G}[cr][cr][0.8][0]{$0.1$}
\psfrag{H}[cr][cr][0.8][0]{$0.2$}
\psfrag{I}[cr][cr][0.8][0]{$0.3$}
\psfrag{J}[cr][cr][0.8][0]{$0.4$}
\psfrag{K}[cr][cr][0.8][0]{$0.5$}
\psfrag{L}[cr][cr][0.8][0]{$0.6$}
\psfrag{M}[cr][cr][0.8][0]{$0.7$}
\psfrag{N}[cr][cr][0.8][0]{$0.8$}
\psfrag{O}[cr][cr][0.8][0]{$0.9$}

\psfrag{Q}[cr][cr][0.8][0]{$\begin{array}{c}\uparrow
\\R^2(t)\end{array}$}
\psfrag{P}[ct][ct][0.8][0]{$t\rightarrow$}

    \begin{tabular}{cc}

      \includegraphics[width=5.5cm]{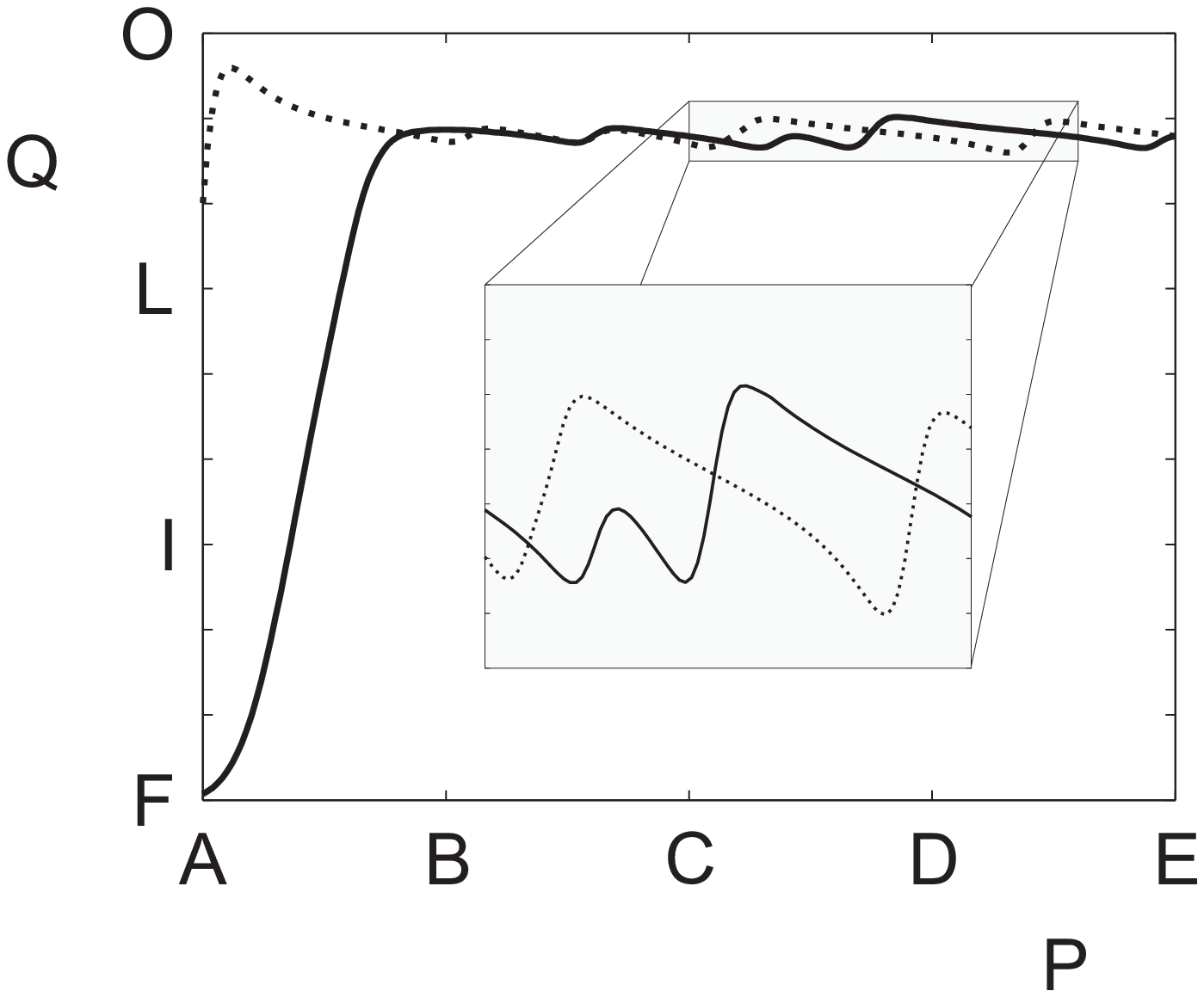}
      & \hspace{0.3cm}
      \includegraphics[width=5.5cm]{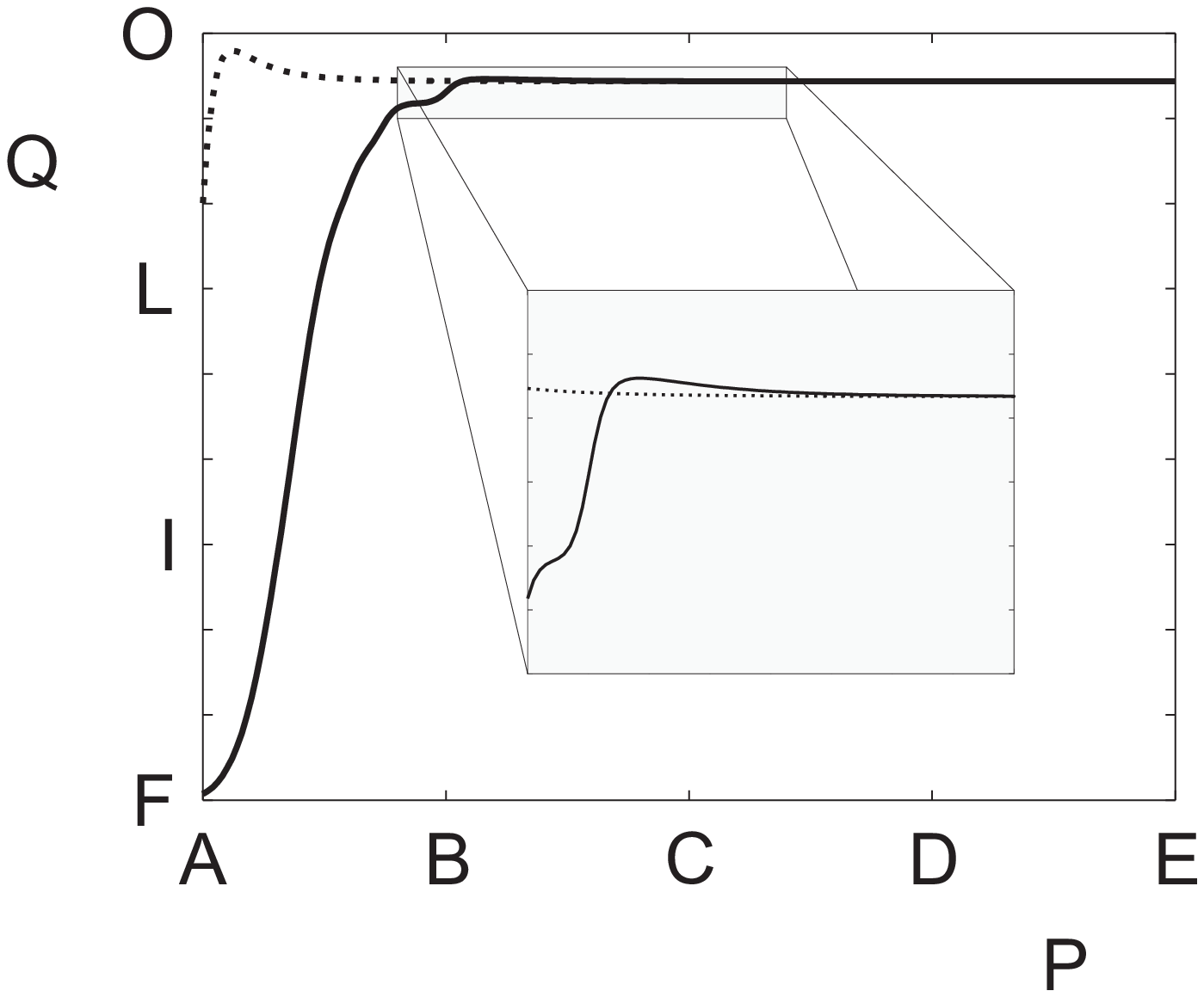}\\
      (a) $k=2.4$ & (b) $k=2.65$

    \end{tabular}

      \caption{Case $N=200$: time evolution of the magnitude squared of the order parameter, $R^2(t)$, for two different initial condition (indicated by a dashed and solid line respectively), and two values of $k$. In the left panel, the value of $k$ ($2.4$) is (well) below the known upper bound on $k_\text{c}$ ($2.6145$)and the system does not converge to a fixed point; in the right panel the value of $k$ ($2.65$) is slightly above the known upper bound on $k_\text{c}$ and the system converges to a fixed point, as expected.}
      \label{fig:Example}
 \end{figure}

 \begin{figure}[t]
\centering \psfrag{A}[ct][ct][0.8][0]{$0$}
\psfrag{B}[ct][ct][0.8][0]{$5$}
\psfrag{C}[ct][ct][0.8][0]{$10$}
\psfrag{D}[ct][ct][0.8][0]{$15$}
\psfrag{E}[ct][ct][0.8][0]{$20$}

\psfrag{F}[cr][cr][0.8][0]{$0$}
\psfrag{I}[cr][cr][0.8][0]{$0.5$}
\psfrag{O}[cr][cr][0.8][0]{$1$}

\psfrag{Q}[cr][cr][0.8][0]{$\begin{array}{c}\uparrow
\\R^2(t)\end{array}$}
\psfrag{P}[ct][ct][0.8][0]{$t\rightarrow$}

    \begin{tabular}{cc}

      \includegraphics[width=5.5cm]{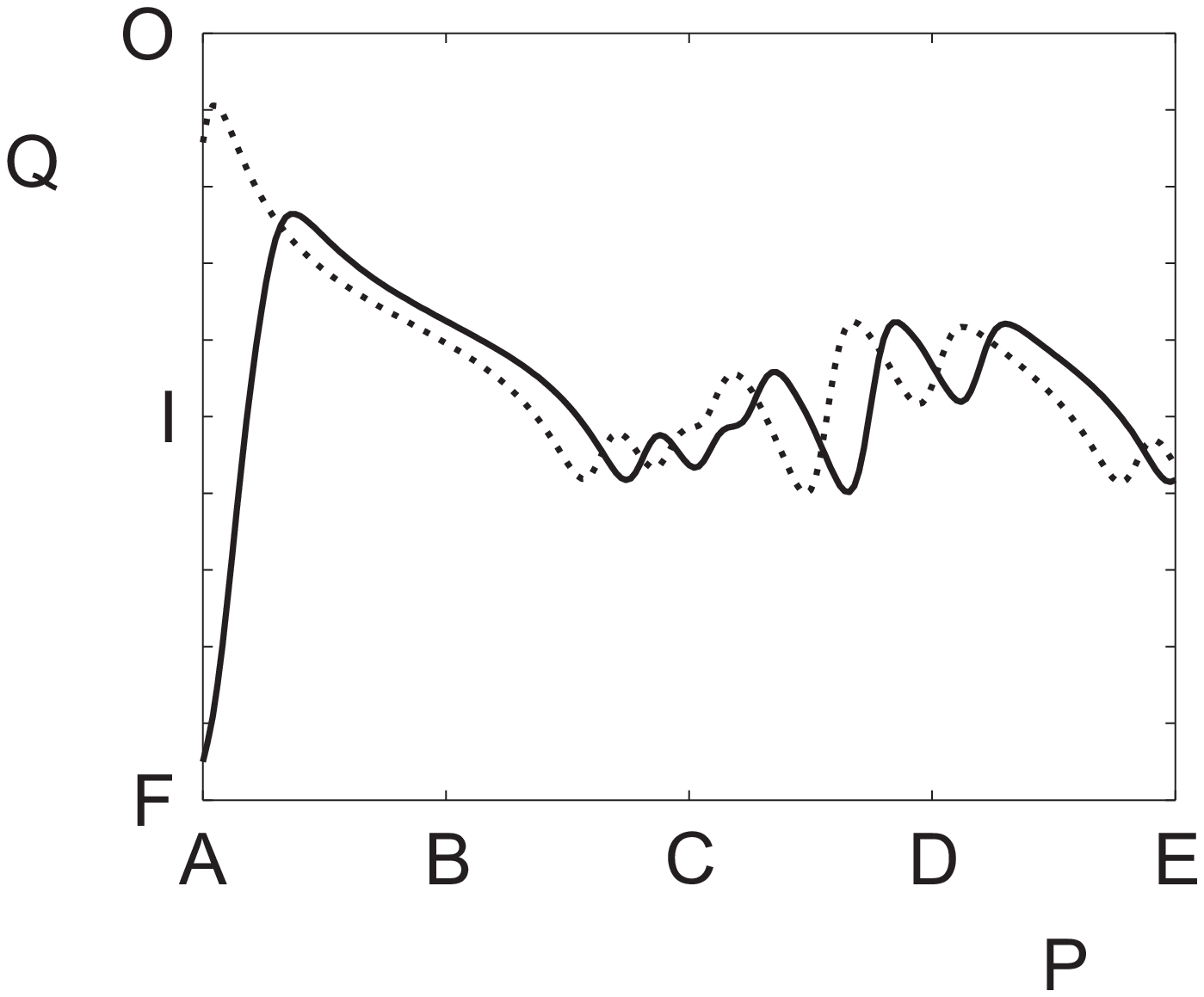}
      & \hspace{0.3cm}
      \includegraphics[width=5.5cm]{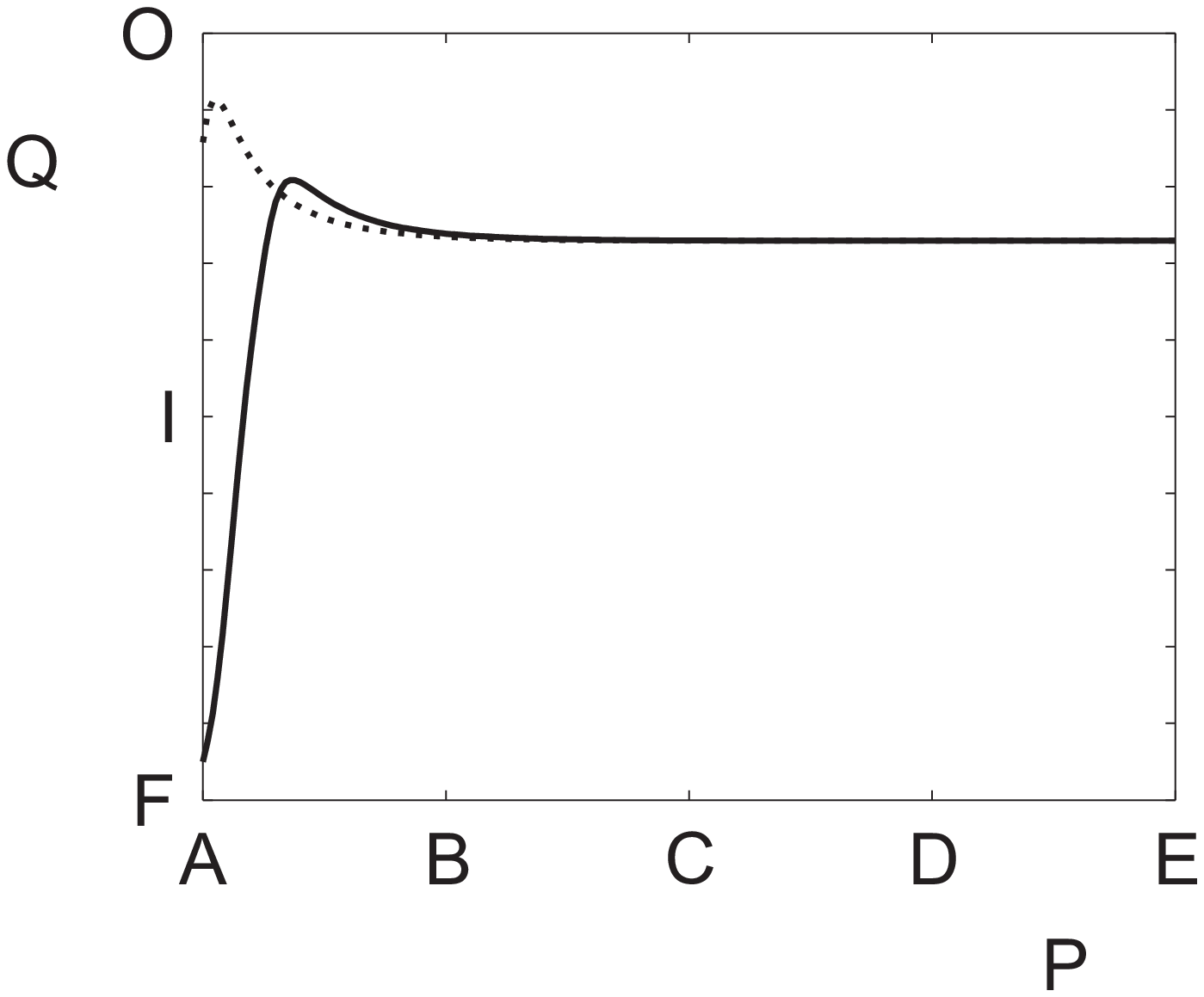}\\
      (a) $k=2.1$ & (b) $k=2.3$

    \end{tabular}

      \caption{Case $N=20$: time evolution of the magnitude squared of the order parameter, $R^2(t)$, for two different initial condition (indicated by a dashed and solid line respectively), and two values of $k$. In the left panel, the value of $k$ ($2.1$) is slightly below the known upper bound on $k_\text{c}$ ($2.2281$)and the system does not converge to a fixed point; in the right panel the value of $k$ ($2.3$) is slightly above the known upper bound on $k_\text{c}$ and the system converges to a fixed point, as expected.}
      \label{fig:Example}
 \end{figure}

\section{Conclusion}We derived necessary and sufficient conditions for the
existence of fixed points in a finite system of coupled
oscillators. In particular, we derived an easy sufficient
condition in terms of the individual oscillator frequencies
(Corollary~\ref{cor:sufficient}), which we used to compute
an upper bound on the critical coupling
(Corollary~\ref{cor:upperbound}). We showed that when no
prior knowledge of the distribution of frequencies is
available, we can still bound the critical coupling in
terms of the infinity norm of the frequencies with their
mean removed (Corollary~\ref{cor:tightbounds}). These
bounds were shown to be the tightest possible, in the sense
that we can find realizations of the intrinsic frequencies
for which the upper bound is attained, and others for which
the critical coupling is arbitrarily close to the lower
bound. Finally, we proposed an efficient algorithm
(Algorithm~\ref{alg:algorithm}) for computing the critical
coupling to within arbitrary bounds in a finite number of
steps. In future work we shall seek to extend the present
analysis to complex networks of arbitrary topology, and
investigate more closely the impact of the shape of the
distribution of intrinsic frequencies on the value of the
critical coupling. We shall also consider the important
question of stability, and present analytical results for
the limit case when the number of oscillators tends to
infinity.

\end{document}